\documentclass[12pt]{article}
\usepackage{amsmath,amssymb}
\usepackage{theorem}
\usepackage{youngtab}
\usepackage{mathdots}
\usepackage{graphicx}
\usepackage{xypic}
\usepackage{epic,eepic}
\usepackage[only,longmapsfrom,trianglelefteqslant]{stmaryrd}
\usepackage{hyperref}
\usepackage{doi}

\newtheorem{theorem}{Theorem}[section]
\newtheorem{proposition}[theorem]{Proposition}
\newtheorem{lemma}[theorem]{Lemma}
\newtheorem{pfpf}{{\it Proof.}}

\newenvironment{prf}{\begin{pfpf}\rm}{\hspace*{\fill}{$\square$}\end{pfpf}}
\theorembodyfont{\rmfamily}\newtheorem{remark}[theorem]{Remark}
\theorembodyfont{\rmfamily}\newtheorem{example}{Example}
\theorembodyfont{\rmfamily}\newtheorem{definition}[theorem]{Definition}

\newcommand{\InjHull}{\operatorname{E}}
\newcommand{\nontrivloops}{\Omega}

\numberwithin{equation}{section}

\begin{document}

\title{\textbf{Young's lattice and dihedral symmetries revisited:
{M}\"obius strips \& metric geometry}\\[9mm]}

\author{\Large Ruedi Suter\\[3mm]
{\textit{Department of Mathematics}}\\
{\textit{ETH Zurich}}\\
{\textit{Raemistrasse 101, 8092 Zurich, Switzerland}}\\[2mm]
e-mail: \texttt{suter@math.ethz.ch}}

\maketitle

\begin{abstract}
A cascade of dihedral symmetries is hidden in Young's lattice of integer
partitions. In fact, for each $N\in\mathbb Z_{\geqslant3}$ the Hasse graph of
the subposet $\mathbb Y_N$ consisting of the partitions with maximal hook
length strictly less than $N$ has the dihedral group of order $2N$ as its
symmetry group. Here a new interpretation of those Hasse graphs is presented,
namely as the $1$-skeleta of the injective hulls of certain finite metric spaces.
\end{abstract}

\section{Introduction and some history}
For each positive integer $N$ let $\mathbb Y_N$ denote the set of those
integer partitions whose maximal hook lengths are strictly less than $N$.
Recall that a partition $\lambda$ covers a partition $\mu$ in Young's
lattice means that the Young diagram of $\mu$ is got by removing an inner
corner box from the Young diagram of $\lambda$. Consider the Hasse
diagram of the subposet of Young's lattice on $\mathbb Y_N$. Let us call
$\operatorname{Hasse}(\mathbb Y_N)$ the underlying undirected abstract graph. 

\begin{example}[Hasse diagram of $\mathbb Y_5$ and its undirected Hasse graph]
\mbox{}
\begin{center}
\mbox{}{\tiny\Yboxdim1mm
\xymatrix@R=2mm@C=2mm{&&&{\yng(2,2,2)}&&
{\yng(3,3)}\\
&&&{\yng(1,2,2)}\ar@{-}[u]&&{\yng(2,3)}\ar@{-}[u]\\
{\yng(1,1,1,1)}&&{\yng(1,1,2)}\ar@{-}[ur]&&{\yng(2,2)}\ar@{-}[ul]\ar@{-}[ur]
&&{\yng(1,3)}\ar@{-}[ul]&&{\yng(4)}\\
&&{\yng(1,1,1)}\ar@{-}[ull]\ar@{-}[u]&&{\yng(1,2)}\ar@{-}[ull]
\ar@{-}[u]\ar@{-}[urr]&&{\yng(3)}\ar@{-}[u]\ar@{-}[urr]\\
&&&{\yng(1,1)}\ar@{-}[ul]\ar@{-}[ur]
&&{\yng(2)}\ar@{-}[ul]\ar@{-}[ur]\\
&&&&{\yng(1)}\ar@{-}[ul]\ar@{-}[ur]\\\\
&&&&{\mbox{}}\ar@{-}[uu]}
\qquad
\qquad
}
\setlength{\unitlength}{0.00088in}
\begin{picture}(2000,0)(-1000,800)
\put(0,0){\circle*{60}}
\put(0,-1000){\circle*{60}}
\put(951,-309){\circle*{60}}
\put(588,809){\circle*{60}}
\put(-588,809){\circle*{60}}
\put(-951,-309){\circle*{60}}
\put(0,-724){\circle*{60}}
\put(688,-224){\circle*{60}}
\put(425,585){\circle*{60}}
\put(-425,585){\circle*{60}}
\put(-688,-224){\circle*{60}}
\put(263,-362){\circle*{60}}
\put(425,138){\circle*{60}}
\put(0,447){\circle*{60}}
\put(-425,138){\circle*{60}}
\put(-263,-362){\circle*{60}}
\thicklines
\path(0,-1000)(0,-724)
\path(951,-309)(688,-224)
\path(588,809)(425,585)
\path(-588,809)(-425,585)
\path(-951,-309)(-688,-224)
\path(0,-724)(263,-362)(688,-224)(425,138)(425,585)(0,447)(-425,585)
(-425,138)(-688,-224)(-263,-362)(0,-724)
\path(0,0)(263,-362)
\path(0,0)(425,138)
\path(0,0)(0,447)
\path(0,0)(-425,138)
\path(0,0)(-263,-362)
\end{picture}
\end{center}
This graph has a $5$-fold (cyclic) symmetry, and its full symmetry group is
a dihedral group of order $10$.
\end{example}

The phenomenon generalizes to all $N$.

\begin{theorem}[see \mbox{\cite[Theorem~2.1]{Su1}}]
For each positive integer $N$ the cyclic group of order $N$ acts faithfully
on $\operatorname{Hasse}(\mathbb Y_N)$.
\end{theorem}

\begin{remark}
Mapping a partition to its dual (also known as conjugate or transpose)
partition restricts to an involution $\mathbb Y_N\to\mathbb Y_N$, which is
different from the identity if $N\geqslant3$. Together with the cyclic symmetry of
order $N$ it generates the full symmetry group $\operatorname{Aut}(
\operatorname{Hasse}(\mathbb Y_N))$, which is a dihedral group of order $2N$ if
$N\geqslant3$.
\end{remark}

The cyclic symmetries emerged as a byproduct of the work on abelian ideals in
a Borel subalgebra $\mathfrak b$ of a finite-dimensional complex simple Lie
algebra $\mathfrak g$ (see \cite{Su2}), namely when $\mathfrak g$ has type
$\mathsf A_{N-1}$.

Here is a quick recapitulation. Let us use the convention that the nilradical
of $\mathfrak b$ is the sum of the root spaces $\mathfrak g_\varphi$
where $\varphi$ runs over the \emph{positive} roots. Then each abelian ideal
$\mathfrak a\trianglelefteqslant\mathfrak b$ is a sum of root spaces
$\mathfrak a=\bigoplus_{\varphi\in\Psi}\mathfrak g_\varphi$ for a certain subset
$\Psi$ of the set of positive roots. If $\rho$ denotes half the sum of the
positive roots, then $\rho+\sum_{\varphi\in\Psi}\varphi$ is the unique integral
weight in the interior of an affine Weyl group translate of the fundamental
alcove (for instance $\rho+\theta$, where $\theta$ is the highest root, is
such an integral weight in the interior of an alcove, which is adjacent
to the fundamental alcove). The union of all those alcoves that one gets if
$\mathfrak a$ runs over all abelian ideals in $\mathfrak b$, is exactly the
fundamental alcove dilated by a factor of two.

In the special case of $\mathfrak{sl}_N(\mathbb C)$ the $N$-fold cyclic
symmetry of the fundamental alcove, which is reflected in the fact that
the affine Dynkin diagram is an $N$-cycle, yields the cyclic symmetry
acting on $\mathbb Y_N$ after identifying elements of $\mathbb Y_N$ in an
evident way with abelian ideals and hence with those alcoves that lie in
(and tessellate) the fundamental alcove dilated by a factor of two.

Since integer partitions and Young's lattice are so basic mathematical
objects, it was natural to look for an independent approach to the
cyclic symmetries. This was realized in \cite{Su1}, and I
reported on it in a talk ``A surprising result about Young's lattice'' in
the ETH Zurich Algebra-Topology Seminar in October 2001. In October
2007, I mentioned the cyclic symmetries parenthetically in a talk
``Some insights from cluster categories'' at INdAM in Rome at a
workshop organized by Paolo Papi and Eric Sommers.
For recent developments see \cite{BZ,TW}.

On February 29, 2012, Urs Lang gave a talk ``Injective hulls of metric
spaces: old and new'' in the ETH Zurich Geometry Seminar. In the introduction
he showed the picture \cite[top of p.~338]{Dr}, which roughly looks as follows

\enlargethispage*{2mm}
\setlength{\unitlength}{0.00048in}
\begin{center}
\begin{picture}(6344,3509)(0,-10)
\texture{55888888 88555555 5522a222 a2555555 55888888 88555555 552a2a2a 2a555555 
	55888888 88555555 55a222a2 22555555 55888888 88555555 552a2a2a 2a555555 
	55888888 88555555 5522a222 a2555555 55888888 88555555 552a2a2a 2a555555 
	55888888 88555555 55a222a2 22555555 55888888 88555555 552a2a2a 2a555555 }
\shade\path(3022,2572)(4522,2722)(5122,1222)
	(3622,1072)(3022,2572)
\path(3022,2572)(4522,2722)(5122,1222)
	(3622,1072)(3022,2572)
\texture{55888888 88555555 5522a222 a2555555 55888888 88555555 552a2a2a 2a555555 
	55888888 88555555 55a222a2 22555555 55888888 88555555 552a2a2a 2a555555 
	55888888 88555555 5522a222 a2555555 55888888 88555555 552a2a2a 2a555555 
	55888888 88555555 55a222a2 22555555 55888888 88555555 552a2a2a 2a555555 }
\shade\path(1822,2122)(3022,2572)(3622,1072)
	(2422,622)(1822,2122)
\path(1822,2122)(3022,2572)(3622,1072)
	(2422,622)(1822,2122)
\texture{55888888 88555555 5522a222 a2555555 55888888 88555555 552a2a2a 2a555555 
	55888888 88555555 55a222a2 22555555 55888888 88555555 552a2a2a 2a555555 
	55888888 88555555 5522a222 a2555555 55888888 88555555 552a2a2a 2a555555 
	55888888 88555555 55a222a2 22555555 55888888 88555555 552a2a2a 2a555555 }
\shade\path(3322,2872)(4822,3022)(4522,2722)
	(3022,2572)(3322,2872)
\path(3322,2872)(4822,3022)(4522,2722)
	(3022,2572)(3322,2872)
\texture{55888888 88555555 5522a222 a2555555 55888888 88555555 552a2a2a 2a555555 
	55888888 88555555 55a222a2 22555555 55888888 88555555 552a2a2a 2a555555 
	55888888 88555555 5522a222 a2555555 55888888 88555555 552a2a2a 2a555555 
	55888888 88555555 55a222a2 22555555 55888888 88555555 552a2a2a 2a555555 }
\shade\path(2722,3172)(3322,2872)(3022,2572)
	(2422,2872)(2722,3172)
\path(2722,3172)(3322,2872)(3022,2572)
	(2422,2872)(2722,3172)
\texture{55888888 88555555 5522a222 a2555555 55888888 88555555 552a2a2a 2a555555 
	55888888 88555555 55a222a2 22555555 55888888 88555555 552a2a2a 2a555555 
	55888888 88555555 5522a222 a2555555 55888888 88555555 552a2a2a 2a555555 
	55888888 88555555 55a222a2 22555555 55888888 88555555 552a2a2a 2a555555 }
\shade\path(2422,2872)(3022,2572)(1822,2122)
	(1222,2422)(2422,2872)
\path(2422,2872)(3022,2572)(1822,2122)
	(1222,2422)(2422,2872)
\thicklines
\path(2422,2872)(2722,3172)
\path(3022,2572)(3322,2872)
\path(3622,1072)(2422,622)
\path(3022,2572)(1822,2122)
\path(1222,2422)(1822,2122)
\path(1822,2122)(2422,622)
\path(3022,2572)(4522,2722)
\path(3322,2872)(4822,3022)
\path(3622,1072)(5122,1222)
\path(4522,2722)(4822,3022)
\path(2722,3172)(3322,2872)
\path(2722,3172)(2722,3472)
\path(1222,2422)(22,3472)
\path(4822,3022)(5572,3472)
\path(2422,622)(2122,22)
\path(5122,1222)(6322,622)
\path(1222,2422)(2422,2872)
\path(4522,2722)(5122,1222)
\path(2422,2872)(3022,2572)
\path(3022,2572)(3622,1072)
\end{picture}
\end{center}
and shows the injective hull (of one of three possible types) of a generic
metric space with five points.

A drawing for the injective hull of a generic metric space with four points
shows a rectangle with four ``antennas'' attached, one at each vertex;
for three points one gets a $\mathsf Y$-shaped tree;
for two points the injective hull is a line segment;
and the injective hull of a one-point space has one point.
Hence for $N\in\{1,2,3,4,5\}$ we recognize the geometric realizations of
the graphs $\operatorname{Hasse}(\mathbb Y_N)$ as the $1$-skeleta of the
injective hulls of certain $N$-point metric spaces.

What about the $N=6$ situation?
From \cite{Su1} the graph $\operatorname{Hasse}(\mathbb Y_6)$
on $32$~vertices and with $48$~edges looks as in the following picture.
\setlength{\unitlength}{0.00044in}
\begin{center}
\begin{picture}(6396,7373)(0,-10)
\put(6316,5479){\circle*{110}}
\put(6316,1879){\circle*{110}}
\put(3198,7279){\circle*{110}}
\put(80,5479){\circle*{110}}
\put(80,1879){\circle*{110}}
\put(3198,79){\circle*{110}}
\put(3198,6079){\circle*{110}}
\put(5276,4879){\circle*{110}}
\put(5276,2479){\circle*{110}}
\put(3198,1279){\circle*{110}}
\put(1120,4879){\circle*{110}}
\put(1120,2479){\circle*{110}}
\put(3198,4879){\circle*{110}}
\put(4237,4279){\circle*{110}}
\put(4237,3079){\circle*{110}}
\put(3198,2479){\circle*{110}}
\put(2159,3079){\circle*{110}}
\put(2159,4279){\circle*{110}}
\put(3198,3499){\circle*{110}}
\put(3198,3859){\circle*{110}}
\put(3718,5779){\circle*{110}}
\put(3718,1579){\circle*{110}}
\put(2678,5779){\circle*{110}}
\put(2678,1579){\circle*{110}}
\put(4756,5179){\circle*{110}}
\put(4756,2179){\circle*{110}}
\put(1640,5179){\circle*{110}}
\put(1640,2179){\circle*{110}}
\put(5276,4279){\circle*{110}}
\put(5276,3079){\circle*{110}}
\put(1120,4279){\circle*{110}}
\put(1120,3079){\circle*{110}}
\thicklines
\path(3198,4879)(2159,4279)(2159,3079)
	(3198,2479)(4237,3079)(4237,4279)(3198,4879)
\path(3198,6079)(3198,7279)
\path(5276,2479)(6316,1879)
\path(3198,1279)(3198,79)
\path(3198,6079)(1120,4879)(1120,2479)
	(3198,1279)(5276,2479)(5276,4879)(3198,6079)
\path(1120,4879)(80,5479)
\path(1120,2479)(80,1879)
\path(5276,4879)(6316,5479)
\path(3198,2479)(3198,3499)
\path(2159,3079)(3198,3859)(4237,3079)
\path(3198,3859)(3198,4879)
\path(4237,4279)(5276,4279)
\path(4237,3079)(5276,3079)
\path(2159,4279)(1120,4279)
\path(2159,3079)(1120,3079)
\path(4237,3079)(4756,2179)
\path(3198,2479)(2678,1579)
\path(2159,3079)(1640,2179)
\path(2159,4279)(3198,3499)(4237,4279)
\path(3198,4879)(3718,5779)
\path(3198,2479)(3718,1579)
\path(3198,4879)(2678,5779)
\path(4237,4279)(4756,5179)
\path(2159,4279)(1640,5179)
\end{picture}
\end{center}
This graph is indeed isomorphic to the $1$-skeleton of the
injective hull of the six-point metric space visualized in
\cite[Fig.~1 on p.~570]{HKM} or in \cite[Fig.~2 on p.~176]{HJ}.

In general, $\operatorname{Hasse}(\mathbb Y_N)$ can be geometrically realized
as the $1$-skeleton of the injective hull of an $N$-point metric space
where the $N$ points form an orbit under an isometric cyclic action.

\begin{theorem}\label{SummaryThm}
Let $N\in\mathbb Z_{\geqslant2}$. Consider a geometric realization of the graph
$\operatorname{Hasse}(\mathbb Y_N)$ as a metric space in which all edges have
length $1$ and let $X_N$ be its boundary, that is, $X_N$ is the $N$-point
subspace consisting of the $N$ pending vertices (the tips of the ``antennas'')
corresponding to the rectangular partitions $(j^{N-j})$ (for $j=1,\dots,N-1$)
together with the empty partition.
Then the $1$-skeleton of the injective hull \,$\InjHull(X_N)$ is a geometric
realization of the graph $\operatorname{Hasse}(\mathbb Y_N)$.
\end{theorem}
For a proof see Remark~\ref{ProofSummaryThm}.

\section{M\"obius strips and Young diagrams}\label{geography}
Let $N\in\mathbb Z_{\geqslant2}$. We consider the discrete circle with
$N$ sites $X=X_N=\{0,1,\dots,N-1\}$. Two sites $j,k\in X$ are neighbouring
if $|k-j|=1$ or $|k-j|=N-1$. Expressed more uniformly, the two sites $j$ and
$j+1$ are neighbouring for $j=0,\dots,N-1$, where $N$ is identified with $0$.

\begin{example}
The circle for $N=4$ with its $4$ sites.
\setlength{\unitlength}{12mm}
\begin{center}
\begin{picture}(12,0.7)(0,-0.3)
\put(0,0){\line(1,0){5}}
\put(7,0){\line(1,0){5}}
\multiput(1,-0.1)(1,0){5}{\line(0,1){0.2}}
\multiput(8,-0.1)(1,0){5}{\line(0,1){0.2}}
\thicklines
\multiput(0,-0.1)(4,0){2}{\line(0,1){0.2}}
\multiput(7,-0.1)(4,0){2}{\line(0,1){0.2}}
\put(0.5,0.07){\makebox(0,0)[b]{\footnotesize$0$}}
\put(1.5,0.07){\makebox(0,0)[b]{\footnotesize$1$}}
\put(2.5,0.07){\makebox(0,0)[b]{\footnotesize$2$}}
\put(3.5,0.07){\makebox(0,0)[b]{\footnotesize$3$}}
\put(4.5,0.07){\makebox(0,0)[b]{\footnotesize$4$}}
\put(6,0){\makebox(0,0){$=$}}
\put(7.5,0.07){\makebox(0,0)[b]{\footnotesize$0$}}
\put(8.5,0.07){\makebox(0,0)[b]{\footnotesize$1$}}
\put(9.5,0.07){\makebox(0,0)[b]{\footnotesize$2$}}
\put(10.5,0.07){\makebox(0,0)[b]{\footnotesize$3$}}
\put(11.5,0.07){\makebox(0,0)[b]{\footnotesize$0$}}
\end{picture}
\end{center}
\end{example}

Next we consider a discrete M\"obius strip $\mathfrak X=\mathfrak X_N$ with
$X$ embedded as its boundary. Namely,
$\mathfrak X=X\times X\!\bigm/\!\{(j,k)\sim(k,j)\mid{j,k}\in X\}$.
The M\"obius strip $\mathfrak X$ has
$\frac12N(N+1)$ sites represented by $(j,k)$ ($0\leqslant j
\leqslant k\leqslant N-1$), and we may extend this to $0\leqslant j
\leqslant k\leqslant N$ by identifying $N$ with $0$.
The embedding $X\hookrightarrow\mathfrak X$ is $j\mapsto(j,j)$.
By convention and slightly abusing the notation, we write
$(j,k)\in\mathfrak X$ for the \emph{class} of the pair $(j,k)$ (i.\,e., its
site in $\mathfrak X$) and hence $(j,k)=(k,j)\in\mathfrak X$. 
Some definitions below (for instance, in Lemma~\ref{Lregular}) use the
independence under the exchange $j\leftrightarrow k$ for their well-definedness,
something that we keep in mind with the tacit understanding. 

\begin{example}
The M\"obius strip for $N=4$ with its $10$ sites. 
\setlength{\unitlength}{0.00040in}
\begin{center}
\begin{picture}(6034,3649)(0,-10)
\thicklines
\path(22,3012)(622,3612)(1222,3012)
	(1822,3612)(2422,3012)(3022,3612)
	(3622,3012)(4222,3612)(4822,3012)
	(2422,612)(22,3012)
\thinlines
\path(2422,612)(3022,12)(5422,2412)(4822,3012)
\path(1222,3012)(3622,612)
\path(2422,3012)(4222,1212)
\path(3622,3012)(4822,1812)
\path(4822,3012)(5422,3612)(6022,3012)(5422,2412)
\path(3622,3012)(1822,1212)
\path(2422,3012)(1222,1812)
\path(1222,3012)(622,2412)
\put(622,3012){\makebox(0,0){\footnotesize$(0,0)$}}
\put(1822,3012){\makebox(0,0){\footnotesize$(1,1)$}}
\put(3022,3012){\makebox(0,0){\footnotesize$(2,2)$}}
\put(4222,3012){\makebox(0,0){\footnotesize$(3,3)$}}
\put(1222,2412){\makebox(0,0){\footnotesize$(0,1)$}}
\put(2422,2412){\makebox(0,0){\footnotesize$(1,2)$}}
\put(3622,2412){\makebox(0,0){\footnotesize$(2,3)$}}
\put(1822,1812){\makebox(0,0){\footnotesize$(0,2)$}}
\put(3022,1812){\makebox(0,0){\footnotesize$(1,3)$}}
\put(2422,1212){\makebox(0,0){\footnotesize$(0,3)$}}
\put(3022,612){\makebox(0,0){\footnotesize$(0,4)$}}
\put(3622,1212){\makebox(0,0){\footnotesize$(1,4)$}}
\put(4222,1812){\makebox(0,0){\footnotesize$(2,4)$}}
\put(4822,2412){\makebox(0,0){\footnotesize$(3,4)$}}
\put(5422,3012){\makebox(0,0){\footnotesize$(4,4)$}}
\end{picture}\makebox[22mm][c]{\raisebox{17mm}[0pt][0pt]{$=$}}
\begin{picture}(6034,3649)(0,-10)
\thicklines
\path(22,3012)(622,3612)(1222,3012)
	(1822,3612)(2422,3012)(3022,3612)
	(3622,3012)(4222,3612)(4822,3012)
	(2422,612)(22,3012)
\thinlines
\path(2422,612)(3022,12)(5422,2412)(4822,3012)
\path(1222,3012)(3622,612)
\path(2422,3012)(4222,1212)
\path(3622,3012)(4822,1812)
\path(4822,3012)(5422,3612)(6022,3012)(5422,2412)
\path(3622,3012)(1822,1212)
\path(2422,3012)(1222,1812)
\path(1222,3012)(622,2412)
\put(622,3012){\makebox(0,0){\footnotesize$(0,0)$}}
\put(1822,3012){\makebox(0,0){\footnotesize$(1,1)$}}
\put(3022,3012){\makebox(0,0){\footnotesize$(2,2)$}}
\put(4222,3012){\makebox(0,0){\footnotesize$(3,3)$}}
\put(1222,2412){\makebox(0,0){\footnotesize$(0,1)$}}
\put(2422,2412){\makebox(0,0){\footnotesize$(1,2)$}}
\put(3622,2412){\makebox(0,0){\footnotesize$(2,3)$}}
\put(1822,1812){\makebox(0,0){\footnotesize$(0,2)$}}
\put(3022,1812){\makebox(0,0){\footnotesize$(1,3)$}}
\put(2422,1212){\makebox(0,0){\footnotesize$(0,3)$}}
\put(3022,612){\makebox(0,0){\footnotesize$(0,0)$}}
\put(3622,1212){\makebox(0,0){\footnotesize$(0,1)$}}
\put(4222,1812){\makebox(0,0){\footnotesize$(0,2)$}}
\put(4822,2412){\makebox(0,0){\footnotesize$(0,3)$}}
\put(5422,3012){\makebox(0,0){\footnotesize$(0,0)$}}
\end{picture}
\end{center}
\end{example}

We need to familiarize ourselves with some notations for the geography of
the M\"obius strip $\mathfrak X$. Let us first have a local inspection.
\begin{itemize}
\item[1)] Each pair of adjacent sites in the M\"obius
strip is represented (in at least one way) as
\setlength{\unitlength}{11mm}
\begin{center}
\begin{picture}(7,3)
\multiput(0,2)(1,-1){3}{\line(1,1){1}}
\multiput(0,2)(1,1){2}{\line(1,-1){2}}
\put(1,2){\makebox(0,0){\footnotesize$(j,k)$}}
\put(2,1){\makebox(0,0){\footnotesize$(j,k+1)$}}
\put(3,1.5){\makebox(0,0)[l]{$(0\leqslant j\leqslant k\leqslant N-1)$\quad or}}
\end{picture}
\quad
\begin{picture}(6,3)
\multiput(0,1)(1,-1){2}{\line(1,1){2}}
\multiput(0,1)(1,1){3}{\line(1,-1){1}}
\put(1,1){\makebox(0,0){\footnotesize$(j,k)$}}
\put(2,2){\makebox(0,0){\footnotesize$(j+1,k)$}}
\put(3,1.5){\makebox(0,0)[l]{$(0\leqslant j<k\leqslant N)$}}
\end{picture}
\end{center}
The pairs on the right for $k=N$ are already represented by
the pairs of sites $(0,j)$ and $(0,j+1)$ depicted on the left.
\item[2)] Each triple of sites consisting of a site adjacent to two sites
at the boundary of the M\"obius strip can be represented as
\begin{center}
\begin{picture}(8,3)
\put(0,2){\line(1,1){1}}
\multiput(1,1)(1,-1){2}{\line(1,1){2}}
\multiput(0,2)(1,1){2}{\line(1,-1){2}}
\put(3,3){\line(1,-1){1}}
\put(1,2){\makebox(0,0){\footnotesize$(j,j)$}}
\put(2,1){\makebox(0,0){\footnotesize$(j,j+1)$}}
\put(3,2){\makebox(0,0){\footnotesize$(j\!+\!1,j\!+\!1)$}}
\put(4,1){\makebox(0,0)[l]{$(0\leqslant j\leqslant N-1)$}}
\end{picture}
\end{center}
\newpage
\item[3)] Each $2$\/$\times$\/$2$ square of sites
in the M\"obius strip can be represented as
\begin{center}
\begin{picture}(8,4)
\multiput(0,2)(1,-1){3}{\line(1,1){2}}
\multiput(0,2)(1,1){3}{\line(1,-1){2}}
\put(1,2){\makebox(0,0){\footnotesize$(j,k)$}}
\put(2,1){\makebox(0,0){\footnotesize$(j,k+1)$}}
\put(2,3){\makebox(0,0){\footnotesize$(j+1,k)$}}
\put(3,2){\makebox(0,0){\footnotesize$(j\!+\!1,k\!+\!1)$}}
\put(4,1){\makebox(0,0)[l]{$(0\leqslant j<k\leqslant N-1)$}}
\end{picture}
\end{center}
(For $N=2$ this square has only three different sites.) 
\end{itemize}

As regards global objects in the geography of $\mathfrak X$, it will soon become
evident that homotopically nontrivial loops are important to consider. Let
$(j_0,k_0)$ adjacent to $(j_1,k_1)$ adjacent to $(j_2,k_2)$ adjacent to $\ldots$
adjacent to $(j_{N-1},k_{N-1})$ adjacent to
$(j_N,k_N)=(j_0,k_0)$ be such a loop consisting of $N$ sites (see
(\ref{exampleloop}) for an example).

\begin{definition}\label{shortnontrivialloops}
Let
$$\nontrivloops=\nontrivloops_N=\Biggl\{\mathcal L\subseteq\mathfrak X\Biggm|
\mbox{\begin{minipage}{8cm}$|\mathcal L|=N$
and the sites in $\mathcal L$ realize a homotopically nontrivial
loop in the M\"obius strip\end{minipage}}\Biggr\}.$$
\end{definition}
Note that for $\mathcal L\subseteq\mathfrak X$ with $|\mathcal L|<N$
the sites in $\mathcal L$ cannot realize a homotopically nontrivial
loop in the M\"obius strip.

\begin{remark}\label{exhaustX}
Since for $\mathcal L\in\nontrivloops$ the sites in $\mathcal L$
realize a homotopically nontrivial loop, the following property holds:
for each $\mathcal L\in\nontrivloops$
$$X=\bigcup_{(j,k)\in\mathcal L}\{j,k\}.$$
Lemma~\ref{Lregular} below gives more precise information. 
\end{remark}

\begin{proposition}\label{loopsandpartitions}
There is a bijection
\begin{align*}
\nontrivloops&\longrightarrow\mathbb Y_N=\{\lambda\mid
\mbox{$\lambda$ is a partition with maximal hook length ${}<N$}\}.\\
\mathcal L_\lambda&\longmapsfrom\lambda
\end{align*}
In particular, $|\nontrivloops|=2^{N-1}$. The loop $\mathcal L_\lambda$
will be referred to as the (loop associated with the) \emph{outer rim}
of $\lambda$. 
\end{proposition}
\begin{prf}
Represent the sites of $\mathfrak X$ in a triangular shape with sites
$(j,k)$ ($0\leqslant j\leqslant k\leqslant N$). The three corners
$(0,0)$ ``upper left'', $(0,N)$ ``bottom'', and $(N,N)$ ``upper right''
of this triangular shape all represent the same site $(0,0)\in\mathfrak X$,
and $(0,k)=(k,N)\in\mathfrak X$.

For $\mathcal L\in\nontrivloops$ we obtain the corresponding partition
$\lambda$ such that the sites below the sites contained in $\mathcal L$
(in the triangular shape) are just the boxes of the Young diagram
(drawn in Russian convention) of $\lambda$. So $\mathcal L$ can be
considered as the outer rim of $\lambda$. (If $\lambda$ has maximal
hook length $h<N-1$, we could consider $\lambda$ as having its nonzero
parts and in addition $N-1-h$ zero parts.)

\begin{example}
Here is an example with $N=9$ for the partition $\lambda=(5,3,3,2)$.
\setlength{\unitlength}{0.00040in}
\begin{center}
\begin{picture}(12034,6659)(0,-10)
\thicklines
\shade\path(3622,2422)(4822,3622)(5422,3022)
	(6022,3622)(7222,2422)(8422,3622)
	(9022,3022)(6022,22)(3622,2422)
\path(3622,2422)(4822,3622)(5422,3022)
	(6022,3622)(7222,2422)(8422,3622)
	(9022,3022)(6022,22)(3622,2422)
\whiten\path(3022,3022)(4822,4822)(5422,4222)
	(6022,4822)(7222,3622)(8422,4822)
	(9622,3622)(9022,3022)(8422,3622)
	(7222,2422)(6022,3622)(5422,3022)
	(4822,3622)(3622,2422)(3022,3022)
\path(3022,3022)(4822,4822)(5422,4222)
	(6022,4822)(7222,3622)(8422,4822)
	(9622,3622)(9022,3022)(8422,3622)
	(7222,2422)(6022,3622)(5422,3022)
	(4822,3622)(3622,2422)(3022,3022)
\put(3622,3022){\makebox(0,0){\footnotesize$(0,5)$}}
\put(4222,3622){\makebox(0,0){\footnotesize$(1,5)$}}
\put(4822,4222){\makebox(0,0){\footnotesize$(2,5)$}}
\put(5422,3622){\makebox(0,0){\footnotesize$(2,6)$}}
\put(6022,4222){\makebox(0,0){\footnotesize$(3,6)$}}
\put(6622,3622){\makebox(0,0){\footnotesize$(3,7)$}}
\put(7222,3022){\makebox(0,0){\footnotesize$(3,8)$}}
\put(7822,3622){\makebox(0,0){\footnotesize$(4,8)$}}
\put(8422,4222){\makebox(0,0){\footnotesize$(5,8)$}}
\put(9022,3622){\makebox(0,0){\footnotesize$(5,9)$}}
\thinlines
\path(4822,1222)(9622,6022)
\path(4222,1822)(8422,6022)
\path(3022,3022)(6022,6022)
\path(2422,3622)(4822,6022)
\path(1822,4222)(3622,6022)
\path(1222,4822)(2422,6022)
\path(622,5422)(1222,6022)
\path(1222,6022)(6622,622)
\path(2422,6022)(7222,1222)
\path(3622,6022)(7822,1822)
\path(4822,6022)(8422,2422)
\path(6022,6022)(9022,3022)
\path(8422,6022)(10222,4222)
\path(9622,6022)(10822,4822)
\path(10822,6022)(11422,5422)
\path(5422,622)(6022,22)(12022,6022)
	(11422,6622)(10822,6022)
\thicklines
\path(22,6022)(622,6622)(1222,6022)
	(1822,6622)(2422,6022)(3022,6622)
	(3622,6022)(4222,6622)(4822,6022)
	(5422,6622)(6022,6022)(6622,6622)
	(7222,6022)(7822,6622)(8422,6022)
	(9022,6622)(9622,6022)(10222,6622)
	(10822,6022)(5422,622)(22,6022)
\thinlines
\path(7222,6022)(9622,3622)
\path(3622,2422)(7222,6022)
\end{picture}
\end{center}
The outer rim of $\lambda$ is
\begin{equation}\label{exampleloop}
\mathcal L_\lambda=\bigl\{(0,5),(1,5),(2,5),(2,6),(3,6),(3,7),(3,8),(4,8),
(5,8)\bigr\}\subseteq\mathfrak X.
\end{equation}

Here is an encoding of $\mathcal L_{\lambda}$ in terms of a graph with vertex set
$X=\{0,\dots,N-1\}$: there is an edge between $j$ and $k$ if and only if
$(j,k)\in\mathcal L_{\lambda}$.
\setlength{\unitlength}{0.00015in}
\begin{center}
\begin{picture}(7120,7270)(0,-10)
\thicklines
\path(6698,3560)(5926,1439)(3971,310)
\path(5926,1439)(297,2431)(1748,702)
	(1748,6418)
\path(1748,6418)(297,4689)
\path(1748,6418)(3971,6810)(5926,5681)
\path(3971,6810)(5926,1439)
\put(6998,3560){\makebox(0,0){\footnotesize\tiny$0$}}
\put(6156,5874){\makebox(0,0){\footnotesize\tiny$4$}}
\put(4023,7105){\makebox(0,0){\footnotesize\tiny$8$}}
\put(1598,6678){\makebox(0,0){\footnotesize\tiny$3$}}
\put(15,4791){\makebox(0,0){\footnotesize\tiny$7$}}
\put(15,2329){\makebox(0,0){\footnotesize\tiny$2$}}
\put(1598,442){\makebox(0,0){\footnotesize\tiny$6$}}
\put(4023,15){\makebox(0,0){\footnotesize\tiny$1$}}
\put(6156,1246){\makebox(0,0){\footnotesize\tiny$5$}}
\end{picture}
\end{center}
Note that in this way we get a connected graph
with $N$ vertices and $N$ edges, hence with exactly one cycle (which can be
a graph-theoretic loop). From this graph we can easily read off a bijection
$\rho:X\to\mathcal L_\lambda$ such that if $\rho(l)=(j,k)$, then $l\in\{j,k\}$,
so that Remark~\ref{exhaustX} follows as a corollary.
In the example above there is no choice for $\rho(0)=(0,5)$, $\rho(1)=(1,5)$,
$\rho(7)=(3,7)$, $\rho(4)=(4,8)$; for the vertices that belong to the
cycle $(5,2,6,3,8)$, there are exactly two choices:
$$\begin{array}{c|cccccc}
l&5&2&6&3&8\\\hline
\rho(l)\rule{0pt}{5mm}&(2,5)&(2,6)&(3,6)&(3,8)&(5,8)\\
\multicolumn{1}{c|}{\mbox{or}}\\
\rho(l)&(5,8)&(2,5)&(2,6)&(3,6)&(3,8)\\
\end{array}$$

For general $\mathcal L_\lambda$ there are always two choices for the
bijection $\rho$ except if the associated graph has a graph-theoretic loop,
in which case there is a unique such bijection $\rho$.
Note also that the cycle always has odd length. In fact, the
edges of the cycle correspond to the sites where the loop $\mathcal L_\lambda$
turns, and the number of those sites is the number
of outer corners (if the maximal hook length of $\lambda$ is strictly
less than $N-1$) or the number of outer corners minus $2$ (if the maximal
hook length of $\lambda$ is $N-1$) plus the number of inner corners
of $\lambda$. The result then follows from
$$\mbox{\#(outer corners of $\lambda$)}-\mbox{\#(inner corners of $\lambda$)}=1.$$

In the example the outer corners of $\lambda$ are at $(2,6)$ and $(3,8)$
as well as at $(0,5)$ and $(5,9)$ (as a site in $\mathfrak X$ we have
$(0,5)=(5,9)$, and $\mathcal L_\lambda$ does not turn there);
and the inner corners of $\lambda$ are at $(1,6)$ [${}\leftrightarrow(2,5)$],
$(2,7)$ [${}\leftrightarrow(3,6)$], and $(4,9)$ [${}\leftrightarrow(5,8)$],
where [${}\leftrightarrow(j+1,k)$] marks the corresponding site where
$\mathcal L_\lambda$ turns around the inner corner $(j,k+1)$ of $\lambda$.
\end{example}

Let us complete the proof of Proposition~\ref{loopsandpartitions}.
In general, to enumerate all elements $\mathcal L$ of $\nontrivloops$,
we count separately for $1\leqslant k\leqslant N$ those $\mathcal L$
that contain both $(0,k)$ and $(1,k)$. There are $\binom{N-1}{k-1}$ chains
of length $N$ consisting of successively adjacent sites starting at
$(1,k)$ and ending at $(k,N)$ (all inside the triangular shape).
A summation over $k$ gives a total of $2^{N-1}$ possibilities.
\end{prf}

\begin{lemma}\label{Lregular}
For each site $L=(j,k)\in\mathfrak X$ let $e_L:=e_j+e_k\in\mathbb R^X$
with $e_i(l):=\delta_{i,l}$. Or in other words, the vector
$\bigl(e_L(l)\bigr)_{l\in X}$ has entries $1$ at positions $j$ and $k$ if
$j\neq k$ or has an entry $2$ at position $j$ if $j=k\in X$, and all other
entries $0$. For each $\mathcal L_\lambda\in\nontrivloops$\ consider the
$N$\/$\times$\/$N$ matrix
$$T_\lambda:=\bigl(e_L(l)\bigr)_{L\in\mathcal L_\lambda,l\in X}.$$
Then $T_\lambda$ is regular, in fact, its determinant is $\pm2$.
\end{lemma}
\begin{prf}
We start with $\mathcal L_{()}=\bigl\{(0,k)\bigm|0\leqslant k<N
\bigr\}$. The matrix $T_{()}=\bigl(e_{(0,k)}(l)\bigr)_{k,l\in X}$ is triangular
with its $(0,0)$ diagonal entry $2$ and the other diagonal entries $1$, hence
with determinant $2$. We proceed by induction. If $\lambda\neq()$, then
$\lambda$ covers a partition $\mu$ with outer rim $\mathcal L_\mu$. The
situation can be depicted as follows:
\begin{center}
$\mathcal L_\mu\supseteq{}$\
\mbox{\setlength{\unitlength}{4mm}\begin{picture}(4,1.8)(0,1.6)
\put(0,2){\line(1,1){1}}
\multiput(1,1)(1,-1){2}{\line(1,1){2}}
\multiput(0,2)(1,1){2}{\line(1,-1){2}}
\put(3,3){\line(1,-1){1}}
\put(1,2){\makebox(0,0){\footnotesize$B$}}
\put(2,1){\makebox(0,0){\footnotesize$A$}}
\put(3,2){\makebox(0,0){\footnotesize$C$}}
\end{picture}}
\quad$\leftrightsquigarrow$\quad
\mbox{\setlength{\unitlength}{4mm}\begin{picture}(4,1.8)(0,1.6)
\multiput(0,2)(1,-1){2}{\line(1,1){2}}
\put(3,1){\line(1,1){1}}
\multiput(1,3)(1,1){2}{\line(1,-1){2}}
\put(0,2){\line(1,-1){1}}
\put(1,2){\makebox(0,0){\footnotesize$B$}}
\put(2,3){\makebox(0,0){\footnotesize$D$}}
\put(3,2){\makebox(0,0){\footnotesize$C$}}
\end{picture}}\ ${}\subseteq\mathcal L_\lambda$
\end{center}
\mbox{}\\
where $A=(j,k+1)$, $B=(j,k)$, $C=(j+1,k+1)$, and $D=(j+1,k)$.
The Young diagram of $\lambda$ is got by adding the box at $A$ to the
Young diagram of $\mu$.
We can write $\mathcal L_\lambda=\bigl(\mathcal L_\mu-\{A\}\bigr)\cup\{D\}$.
Since $e_D=-e_A+e_B+e_C$, the determinants of the matrices
$T_\lambda$ and $T_\mu$ are equal up to a sign.
(Using the fact that there are exactly one or two distinguished bijections
$\rho:X\to\mathcal L_\lambda$ (with $l\in\{j,k\}$ if $\rho(l)=(j,k)$) and that
if there are two such bijections, then they differ only by an odd-length
cyclic (hence even) permutation, we could remove the ambiguity in the sign of
the determinant.)
\end{prf}

\begin{remark}\label{cyclicaction}
The cyclic action that is given by translation in the M\"obius strip
$\mathfrak X\to\mathfrak X$,
$(j,k)\mapsto(j+1,k+1)$ induces an action on $\nontrivloops$ and
hence on the set of partitions with maximal hook length ${}<N$.
This corresponds to the diagonal sliding operation described
in \cite{Su1} and will be recalled later (see Theorem~\ref{zerocells}).
\end{remark}

\section{Injective hulls of finite metric spaces}
A metric space $Z$ is called \emph{injective} if every $1$-Lipschitz (i.\,e.,
distance non-increasing) map $f:A\to Z$ from a subspace $A$ of any metric
space $X$ can be extended to a $1$-Lipschitz map $\overline f:X\to Z$.

It was first proved by Isbell \cite{Is} that for every metric space $X$ there
is an injective metric space $E$ such that $X$ embeds isometrically into $E$;
and given such an embedding $e:X\to E$, there is a unique smallest injective
subspace $\InjHull(X) \subseteq E$ containing $e(X)$. Let
$\operatorname{e}=\overline{\operatorname{id}_{\InjHull(X)}}\circ e$ where $
\overline{\operatorname{id}_{\InjHull(X)}}:E\to\InjHull(X)$ extends
$\operatorname{id}_{\InjHull(X)}$ as a $1$-Lipschitz map. 
This space $\InjHull(X)$ (or more properly the isometric embedding
$\operatorname{e}:X\to\InjHull(X)$) is called an \emph{injective hull}
of $X$. If $\operatorname{e}':X\to\InjHull'(X)$ is another injective
hull of $X$, then there is an isometry $\iota:\InjHull(X)\to\InjHull'(X)$
with $\operatorname{e}'=\iota\circ\operatorname{e}$. 

The injective hull $\InjHull(X)$ of a finite metric space $X$ with metric $d$
can be realized as the polyhedral complex that consists of the bounded faces
of the polyhedron
\begin{align}\label{Polyhedron}
\Delta(X)&:=\bigl\{f\in\mathbb R^X\bigm|\forall x,y\in X:
f(x)+f(y)\geqslant d(x,y)\bigr\}
\intertext{and can be shown to be (see \cite[Lemma~1]{Dr2})}
\label{ModelInjHull}
\InjHull(X)&:=\bigl\{f\in\Delta(X)\bigm|\forall x\in X\,\exists y\in X:
f(x)+f(y)=d(x,y)\bigr\}.
\end{align}
The distance between two functions $f,g\in\InjHull(X)\subseteq\mathbb R^X$
is $\Vert f-g\Vert_\infty=\max_{x\in X}\bigl|f(x)-g(x)\bigr|$, and $X$ embeds
into $\InjHull(X)$ by $x\mapsto\operatorname{e}(x)=d(x,\phantom{y})$.
\subsection*{$X=X_N=\{0,1,\dots,N-1\}$}
We consider the metric space $X=X_N=\{0,1,\dots,N-1\}$ with metric
$d(0,j)=j(N-j)$ (for $j=0,\dots,N-1$;
the formula holds also if we insert $j=N$ and identify $N$ with $0$) and
extended cyclically, that is, $d(j,k)=|k-j|(N-|k-j|)$.

Note that if we let $0$ correspond to the empty partition and $j\in X-\{0\}$
to the rectangular partition $(j^{N-j})$, then this is the metric induced from
the Hasse graph (all edges of length $1$) of Young's lattice. In fact,
starting with the Young diagram of the partition $(j^{N-j})$, we successively
remove boxes till we get $(j^{N-k})$ (for $k>j$) and then successively add
boxes till we obtain $(k^{N-k})$. 
\setlength{\unitlength}{0.0005in}
\begin{center}
\begin{picture}(5124,3942)(0,-10)
\path(1018.066,551.360)(912.000,615.000)(975.640,508.934)
\path(912,615)(1512,15)
\path(1405.934,78.640)(1512.000,15.000)(1448.360,121.066)
\path(2775.640,121.066)(2712.000,15.000)(2818.066,78.640)
\path(2712,15)(4212,1515)
\path(4148.360,1408.934)(4212.000,1515.000)(4105.934,1451.360)
\path(312,2415)(1812,3915)(3612,2115)
	(2112,615)(312,2415)
\path(1512,1215)(4212,3915)(4812,3315)
	(2112,615)(1512,1215)
\path(118.066,2051.360)(12.000,2115.000)(75.640,2008.934)
\path(12,2115)(1812,315)
\path(1705.934,378.640)(1812.000,315.000)(1748.360,421.066)
\path(2475.640,421.066)(2412.000,315.000)(2518.066,378.640)
\path(2412,315)(5112,3015)
\path(5048.360,2908.934)(5112.000,3015.000)(5005.934,2951.360)
\put(3612,515){\makebox(0,0)[lb]{$j$}}
\put(712,75){\makebox(0,0){$N-k$}}
\put(312,1215){\makebox(0,0)[t]{$N-j$}}
\put(4612,2115){\makebox(0,0)[lb]{$k$}}
\end{picture}
\end{center}
The total number of required moves is
$$\underbrace{j\bigl((N-j)-(N-k)\bigr)}_{\mbox{\footnotesize remove boxes}}
+\underbrace{(k-j)(N-k)}_{\mbox{\footnotesize add boxes}}=(k-j)\bigl(N-(k-j)\bigr)
=d(j,k).$$

\begin{definition}[Extension from $X$ to $\mathfrak X$: $f\mapsto\tilde f$]
\label{embedfX}
We embed $\mathbb R^X\hookrightarrow\mathbb R^{\mathfrak X}$ by
$$f\mapsto \tilde f:(j,k)\mapsto \tfrac12\bigl(f(j)+f(k)-d(j,k)\bigr).$$
This is of course well-defined and $f(j)=\tilde f(j,j)$.
Note also that the extension commutes with affine combinations, that is,
for $(f_i)_{i\in I}\subseteq\mathbb R^X$ and $(a_i)_{i\in I}\subseteq\mathbb R$
with $\sum_{i\in I}a_i=1$ we have
$$\Bigl(\sum_{i\in I}a_if_i\Bigr)\widetilde{\phantom{\bigl(}}
=\sum_{i\in I}a_i\tilde f_i.$$
This applies in particular to convex combinations, where all the
coefficients $a_i$ are non\-negative and sum up to~$1$.
\end{definition}

Let us rewrite the definitions (\ref{Polyhedron}) and (\ref{ModelInjHull})
as in the following definition.
\begin{definition}[Injective hull of $X$]\label{InjHullX}
\begin{align*}
\Delta(X)&:=\bigl\{f\in\mathbb R^X\bigm|\forall L\in\mathfrak X:
\tilde f(L)\geqslant0\bigr\},\\
\InjHull(X)&:=\bigl\{f\in\Delta(X)\bigm|\forall j\in X\,\exists
L=(j,k)=(k,j)\in\mathfrak X:\tilde f(L)=0\bigr\}.
\end{align*}
\end{definition}

\begin{lemma}\label{fliftrecursion}
Let $f\in\mathbb R^X$ and let $\tilde f\in\mathbb R^{\mathfrak X}$
be its extension to $\mathfrak X$ as in Definition~\ref{embedfX}.

\noindent
For any $2$\/$\times$\/$2$ square of sites
\mbox{\ \setlength{\unitlength}{4mm}\begin{picture}(4,2.4)(0,1.6)
\multiput(0,2)(1,-1){3}{\line(1,1){2}}
\multiput(0,2)(1,1){3}{\line(1,-1){2}}
\put(1,2){\makebox(0,0){\footnotesize$B$}}
\put(2,1){\makebox(0,0){\footnotesize$A$}}
\put(2,3){\makebox(0,0){\footnotesize$D$}}
\put(3,2){\makebox(0,0){\footnotesize$C$}}
\end{picture}}${}\subseteq\mathfrak X$ we have\\
\begin{align*}
\tilde f(A)+\tilde f(D)-\tilde f(B)-\tilde f(C)&=1.
\intertext{For any triple of sites
\mbox{\ \setlength{\unitlength}{4mm}\begin{picture}(4,1.6)(0,1.6)
\put(0,2){\line(1,1){1}}
\multiput(1,1)(1,-1){2}{\line(1,1){2}}
\multiput(0,2)(1,1){2}{\line(1,-1){2}}
\put(3,3){\line(1,-1){1}}
\put(1,2){\makebox(0,0){\footnotesize$B$}}
\put(2,1){\makebox(0,0){\footnotesize$A$}}
\put(3,2){\makebox(0,0){\footnotesize$C$}}
\end{picture}}${}\subseteq\mathfrak X$
with $B$ and $C$ at the boundary we have}
2\tilde f(A)+N-\tilde f(B)-\tilde f(C)&=1.
\end{align*}
\end{lemma}
\begin{prf}
Recall from Section~\ref{geography} how to represent the sites in such triples
and quadruples.
For the triples let $A=(j,j+1)$, $B=(j,j)$, $C=(j+1,j+1)$ and compute
\begin{align*}
\makebox[1cm][l]{$2\tilde f(A)+N-\tilde f(B)-\tilde f(C)$}\\
&=f(j)+f(j+1)-(N-1)+N-f(j)-f(j+1)=1
\intertext{and for the quadruples let $A=(j,k+1)$, $B=(j,k)$,
$C=(j+1,k+1)$, $D=(j+1,k)$ and compute}
\makebox[1cm][l]{$\tilde f(A)+\tilde f(D)-\tilde f(B)-\tilde f(C)$}\\
&=\tfrac12\bigl(f(j)+f(k+1)-(k-j+1)(N-k+j-1)\bigr)\\
&\quad\,{}+\tfrac12\bigl(f(j+1)+f(k)-(k-j-1)(N-k+j+1)\bigr)\\
&\quad\,{}-\tfrac12\bigl(f(j)+f(k)-(k-j)(N-k+j)\bigr)\\
&\quad\,{}-\tfrac12\bigl(f(j+1)+f(k+1)-(k-j)(N-k+j)\bigr)=1.\\[-12mm]
\end{align*}
\end{prf}

\begin{definition}\label{HLdjk}
For a site $L=(j,k)\in\mathfrak X$ consider the affine hyperplane
$$H_L:=\bigl\{f\in\mathbb R^X\bigm|f(j)+f(k)=d(j,k)\bigr\}
=\bigl\{f\in\mathbb R^X\bigm|\tilde f(L)=0\bigr\}$$
with $\tilde f$ the extension of $f$ to $\mathfrak X$ as in
Definition~\ref{embedfX}.\par
For a set of sites $\mathcal L\subseteq\mathfrak X$ consider the
intersection of the affine hyperplanes
$$H_{\mathcal L}:=\bigcap_{L\in\mathcal L}H_L\subseteq\mathbb R^X.$$
\end{definition}

The next proposition will be superseded by Theorem~\ref{zerocells}.

\begin{proposition}[Empty partition]
Let $f_{()}:=d(0,\phantom{k})\in\mathbb R^X$ and let $\mathcal L_{()}
\in\nontrivloops$ be the outer rim (see Proposition~\ref{loopsandpartitions})
of the empty partition $()$. Then
\begin{itemize}
\item $H_{\mathcal L_{()}}=\{f_{()}\}$
\item For $L\in\mathfrak X$, say $L=(j,k)$ with $0\leqslant j\leqslant k<N$
(or even $0\leqslant j\leqslant k\leqslant N$), we have
$\tilde f_{()}(L)=\tilde f_{()}(j,k)=j(N-k)$
\item $f_{()}\in\InjHull(X)$
\end{itemize}
\end{proposition}
\begin{prf}
We have $f_{()}(l)=l(N-l)$ and hence for $0\leqslant j\leqslant k<N$
$$\tilde f_{()}(j,k)=\tfrac12\bigl(j(N-j)+k(N-k)-(k-j)(N-(k-j))\bigr)=j(N-k).$$
Since $\mathcal L_{()}=\bigl\{(0,k)\bigm|0\leqslant k<N\bigr\}$ this implies
in particular
\begin{equation}\label{fvaluesempty}
\begin{cases}
\tilde f_{()}(L)=0&\mbox{if $L\in\mathcal L_{()}$}\\
\tilde f_{()}(L)>0&\mbox{if $L\notin\mathcal L_{()}$}
\end{cases}
\end{equation}
and by the regularity expressed in Lemma~\ref{Lregular} (here we need only the 
base case in its inductive proof) we conclude that $H_{\mathcal L_{()}}=\{f_{()}\}$.

The (in)equalities (\ref{fvaluesempty}) show that
$f_{()}\in\Delta(X)$, and to conclude that $f_{()}\in\InjHull(X)$
just take $y=0$ in (\ref{ModelInjHull}).
\end{prf}

\begin{theorem}\label{zerocells}
Consider the $N$-point metric space $X=\{0,1,\dots,N-1\}$ with
metric $d(j,k)=|k-j|(N-|k-j|)$.
Then the $0$-faces of its injective hull $\InjHull(X)$, realized as in
Definition~\ref{InjHullX}, are precisely
$$H_{\mathcal L_\lambda}=\bigcap_{L\in\mathcal L_\lambda}\bigl\{
f\in\mathbb R^X\bigm|\tilde f(L)=0\bigr\}=:\{f_\lambda\}$$
where $\lambda$ runs through the partitions in $\mathbb Y_N$,
and $\mathcal L_\lambda$ is the outer rim of $\lambda$ (as defined in
Proposition~\ref{loopsandpartitions}).
The solutions $f_\lambda\in\mathbb R^X$ have more explicit descriptions:
\begin{enumerate}
\item $f_\lambda$ is determined by the following recursive procedure:
$f_{()}(l)=l(N-l)$; and if $\lambda$ covers $\mu$ such that
$\mathcal L_\mu\ni(j,k+1)\notin\mathcal L_\lambda$ for some pair $(j,k)$
with $0\leqslant j<k\leqslant N-1$, then
\begin{equation}\label{recursivef}
f_\lambda(l)=f_\mu(l)+\begin{cases}
\phantom{+}1&\mbox{for $l=0,\dots,j$,}\\
-1&\mbox{for $l=j+1,\dots,k$,}\\
\phantom{+}1&\mbox{for $l=k+1,\dots,N-1$.}
\end{cases}
\end{equation}
\item $f_\lambda$ can be described directly via $f_\lambda(0)=|\lambda|$ and by
using the cyclic action. For $\lambda=(\lambda_1,\dots,\lambda_m)\in\mathbb Y_N$
let $\lambda_{m+1}=\dots=\lambda_{N-\lambda_1}=0$
(for the empty partition put $m=1$, $\lambda_1=0$). Then
$$\sigma(\lambda)=(\lambda_2+1,\dots,\lambda_{N-\lambda_1}+1)$$
is the cyclic action alluded to in Remark~\ref{cyclicaction}.
Its inverse $\tau$ (note the intended coincidence with the common notation for
an Auslander-Reiten translation) is then $\tau(\lambda)=\bigl(\sigma(\lambda')
\bigr)'$ and explicitly
$$\tau(\lambda)=(N-m-1,\lambda_1-1,\dots,\lambda_m-1)|_{\scriptstyle
\textup{remove trailing zeros}}$$
(for the empty partition put $m=0$, $\lambda_1=1$).
The formula for $f_\lambda$ is
\begin{equation}\label{fdirectly}
f_\lambda(l)=\bigl|\tau^l(\lambda)\bigr|.
\end{equation}
\end{enumerate}
\end{theorem}
\begin{prf}
Let $\lambda\in\mathbb Y_N$. 
By Lemma~\ref{Lregular} the $N$ affine hyperplanes
$\bigl(H_L\bigr)_{L\in\mathcal L_\lambda}$ (Definition~\ref{HLdjk})
intersect in a single point $f_\lambda$.

For a site $L\in\mathfrak X-\mathcal L_\lambda$ consider the partition
$\alpha_\lambda(L)$ which is defined as follows: its Young diagram has the
shape that is bordered by $\mathcal L_\lambda$ and whose ``bottom box''
is at the position of the site $L$.
For $L\in\mathcal L_\lambda$ let $\alpha_\lambda(L)$ be the empty partition.
\setlength{\unitlength}{0.0007in}
\begin{center}
\begin{picture}(3774,2139)(0,-10)
\path(1512,162)(1662,312)(1812,162)
\texture{55888888 88555555 5522a222 a2555555 55888888 88555555 552a2a2a 2a555555 
	55888888 88555555 55a222a2 22555555 55888888 88555555 552a2a2a 2a555555 
	55888888 88555555 5522a222 a2555555 55888888 88555555 552a2a2a 2a555555 
	55888888 88555555 55a222a2 22555555 55888888 88555555 552a2a2a 2a555555 }
\shade\path(12,1362)(462,1812)(612,1662)
	(912,1962)(1212,1662)(1362,1812)
	(1812,1362)(2112,1662)(2262,1512)
	(2712,1962)(3012,1662)(3462,2112)
	(3762,1812)(3612,1662)(3462,1812)
	(3012,1362)(2712,1662)(2262,1212)
	(2112,1362)(1812,1062)(1362,1512)
	(1212,1362)(912,1662)(612,1362)
	(462,1512)(162,1212)(12,1362)
\path(12,1362)(462,1812)(612,1662)
	(912,1962)(1212,1662)(1362,1812)
	(1812,1362)(2112,1662)(2262,1512)
	(2712,1962)(3012,1662)(3462,2112)
	(3762,1812)(3612,1662)(3462,1812)
	(3012,1362)(2712,1662)(2262,1212)
	(2112,1362)(1812,1062)(1362,1512)
	(1212,1362)(912,1662)(612,1362)
	(462,1512)(162,1212)(12,1362)
\path(312,1362)(1662,12)(3012,1362)
\put(1662,162){\makebox(0,0){\tiny$L$}}
\put(1662,800){\makebox(0,0){\small shape $\alpha_\lambda(L)$}}
\put(-75,1125){\makebox(0,0){$\mathcal L_\lambda$}}
\end{picture}
\qquad
\begin{picture}(3774,2139)(0,-10)
\path(1512,162)(1662,312)(1812,162)
\texture{55888888 88555555 5522a222 a2555555 55888888 88555555 552a2a2a 2a555555 
	55888888 88555555 55a222a2 22555555 55888888 88555555 552a2a2a 2a555555 
	55888888 88555555 5522a222 a2555555 55888888 88555555 552a2a2a 2a555555 
	55888888 88555555 55a222a2 22555555 55888888 88555555 552a2a2a 2a555555 }
\shade\path(12,1362)(462,1812)(612,1662)
	(912,1962)(1212,1662)(1362,1812)
	(1812,1362)(2112,1662)(2262,1512)
	(2712,1962)(3012,1662)(3462,2112)
	(3762,1812)(3612,1662)(3462,1812)
	(3012,1362)(2712,1662)(2262,1212)
	(2112,1362)(1812,1062)(1362,1512)
	(1212,1362)(912,1662)(612,1362)
	(462,1512)(162,1212)(12,1362)
\path(12,1362)(462,1812)(612,1662)
	(912,1962)(1212,1662)(1362,1812)
	(1812,1362)(2112,1662)(2262,1512)
	(2712,1962)(3012,1662)(3462,2112)
	(3762,1812)(3612,1662)(3462,1812)
	(3012,1362)(2712,1662)(2262,1212)
	(2112,1362)(1812,1062)(1362,1512)
	(1212,1362)(912,1662)(612,1362)
	(462,1512)(162,1212)(12,1362)
\path(312,1362)(1662,12)(3012,1362)
\path(612,1362)(1662,312)(2862,1512)
\path(1362,312)(1662,612)(1962,312)
\put(1662,162){\makebox(0,0){\tiny$A$}}
\put(1512,312){\makebox(0,0){\tiny$B$}}
\put(1812,312){\makebox(0,0){\tiny$C$}}
\put(1662,462){\makebox(0,0){\tiny$D$}}
\put(-75,1125){\makebox(0,0){$\mathcal L_\lambda$}}
\end{picture}
\end{center}
 In particular, we recover $\lambda$ as
\begin{equation}\label{recoverlambda}
\lambda=\alpha_\lambda(L)\bigr|_{L=(0,0)=(0,N)}.
\end{equation}
For a quadruple $(A,B,C,D)$ of sites as depicted on the right side in the
illustration above we get of course by considering the areas
\begin{equation}\label{ABCDformula}
|\alpha_\lambda(A)|+|\alpha_\lambda(D)|-|\alpha_\lambda(B)|-|\alpha_\lambda(C)|=1
\end{equation}
so that finally
\begin{equation}\label{numberofsites}
\tilde f_\lambda(L)=|\alpha_\lambda(L)|.
\end{equation}
In fact, for $L\in\mathcal L_\lambda$ both sides in (\ref{numberofsites})
vanish. The formula
(\ref{ABCDformula}) computes the values $|\alpha_\lambda(L)|$ for all
$L\in \mathfrak X$ recursively from the values $|\alpha_\lambda(L)|=0$
for $L\in\mathcal L_\lambda$. On the other hand, by
Lemma~\ref{fliftrecursion}, $\tilde f_\lambda(L)$ satisfies the same
recursive formula. Hence (\ref{numberofsites}) holds true.

\begin{remark}
From (\ref{recoverlambda}) and (\ref{numberofsites})
we get the formula (\ref{fdirectly}) for the case $l=0$, namely,
$f_\lambda(0)=|\lambda|$. More generally,
if $\lambda=(\lambda_1,\dots,\lambda_m)$ with dual partition
$\lambda'=(\lambda_1',\dots,\lambda_{m'}')$, then
\begin{align*}
\tilde f_\lambda(0,k)&=\lambda_{k+1}'+\dots+\lambda_{m'}'&&\mbox{($0\leqslant k
\leqslant m'$),}\\
\tilde f_\lambda(0,N-k)&=\lambda_{k+1}+\dots+\lambda_m&&\mbox{($0\leqslant k
\leqslant m$).}
\end{align*}
\end{remark}

The formula (\ref{fdirectly}) follows from
$\alpha_{\tau(\lambda)}(L)=\alpha_{\lambda}(\sigma_{\mathfrak X}(L))$, where
$\sigma_{\mathfrak X}:\mathfrak X\to\mathfrak X$ is the translation
$(j,k)\mapsto(j+1,k+1)$ in the M\"obius strip,
together with (\ref{recoverlambda}) and (\ref{numberofsites}).

Before we continue with the first part of Theorem~\ref{zerocells}, we use
the geometric interpretation to compare $\tilde f_\lambda$ and $\tilde f_\mu$ if
$\lambda$ covers $\mu$. Their outer rims $\mathcal L_\lambda$
and $\mathcal L_\mu$ satisfy
$$\mathcal L_\lambda=\bigl(\mathcal L_\mu-\{(j,k+1)\}\bigr)\cup\{(j+1,k)\}.$$
The following picture illustrates the situation.

\setlength{\unitlength}{0.00083333in}
\begin{center}
\begin{picture}(4824,5139)(0,-10)
\texture{55888888 88555555 5522a222 a2555555 55888888 88555555 552a2a2a 2a555555 
	55888888 88555555 55a222a2 22555555 55888888 88555555 552a2a2a 2a555555 
	55888888 88555555 5522a222 a2555555 55888888 88555555 552a2a2a 2a555555 
	55888888 88555555 55a222a2 22555555 55888888 88555555 552a2a2a 2a555555 }
\shade\path(1812,2712)(2112,3012)(2412,2712)
	(2112,2412)(1812,2712)
\path(1812,2712)(2112,3012)(2412,2712)
	(2112,2412)(1812,2712)
\shade\path(2412,2712)(2712,3012)(3012,2712)
	(2712,2412)(2412,2712)
\path(2412,2712)(2712,3012)(3012,2712)
	(2712,2412)(2412,2712)
\path(312,4812)(4812,312)(4512,12)
	(3912,612)(3612,312)(3912,12)(4512,612)
\path(312,4812)(612,5112)(1212,4512)
	(1512,4812)(1212,5112)(612,4512)
\path(1212,312)(1512,12)(1812,312)
	(2112,12)(2412,312)(2712,12)
	(3012,312)(3312,12)(3612,312)
\path(1512,4812)(1812,5112)(2112,4812)
	(2412,5112)(2712,4812)(3012,5112)(3312,4812)
\path(312,612)(912,12)(1212,312)
	(912,612)(312,12)(12,312)
	(4512,4812)(4212,5112)(3612,4512)
	(3312,4812)(3612,5112)(4212,4512)
\path(2112,2412)(2412,2112)(2712,2412)
\path(2112,3012)(2412,3312)(2712,3012)
\put(612,4812){\makebox(0,0){\tiny$(j\!\!+\!\!1,\!j\!\!+\!\!1)$}}
\put(4212,4812){\makebox(0,0){\tiny$(k,k)$}}
\put(312,312){\makebox(0,0){\tiny$(k\!\!+\!\!1,\!k\!\!+\!\!1)$}}
\put(4512,312){\makebox(0,0){\tiny$(j,j)$}}
\put(2412,2412){\makebox(0,0){\tiny$(j,k\!+\!1)$}}
\put(2412,3012){\makebox(0,0){\tiny$(j\!+\!1,k)$}}
\put(3000,3012){\makebox(0,0)[l]{$\mathcal L_\lambda\ni(j+1,k)\notin
\mathcal L_\mu$}}
\put(3000,2412){\makebox(0,0)[l]{$\mathcal L_\mu\ni(j,k+1)\notin
\mathcal L_\lambda$}}
\put(2412,4012){\makebox(0,0){$\tilde f_\mu(L)=\tilde f_\lambda(L)+1$}}
\put(2412,3712){\makebox(0,0){for $L$ in this sector}}
\put(2412,1412){\makebox(0,0){$\tilde f_\lambda(L)=\tilde f_\mu(L)+1$}}
\put(2412,1112){\makebox(0,0){for $L$ in this sector}}
\put(-50,3012){\makebox(0,0)[l]{$\tilde f_\lambda(L)=\tilde f_\mu(L)$}}
\put(-50,2712){\makebox(0,0)[l]{for $L$ outside the}}
\put(-50,2412){\makebox(0,0)[l]{sectors (as a site in $\mathfrak X$)}}
\end{picture}
\end{center}
Hence $\tilde f_\lambda(L)-\tilde f_\mu(L)\in\{1,-1,0\}$
depending on the position of the site $L$. 
In particular, this proves the formula (\ref{recursivef}) for the boundary
values.

Let us continue with the proof of the first part of Theorem~\ref{zerocells}. 
From (\ref{numberofsites}) we have $f_\lambda\in\Delta(X)$.
Remark~\ref{exhaustX} say that for each $j\in X$ there exists $k\in X$
such that the site $L:=(j,k)\in\mathcal L_\lambda\subseteq\mathfrak X$, and
hence $\tilde f_\lambda(L)=0$, which shows that $f_\lambda\in\InjHull(X)$.

It is clear that $\mathbb Y_N\ni\lambda\mapsto f_\lambda\in\mathbb R^X$ is an
injective function, for instance because
$\mathcal L_\lambda=\bigl\{L\in\mathfrak X\bigm|\tilde f_\lambda(L)=0\bigr\}$
determines $\lambda$ as in (\ref{recoverlambda}).
Hence
we exhibited $\bigl|\mathbb Y_N\bigr|=2^{N-1}$ $0$-faces in $\InjHull(X)$,
which is the maximal possible number of $0$-faces in the injective hull
of a metric space with $N$ points (the $v=0$ case of the following
Theorem~\ref{HJThm} by Herrmann and Joswig).
Hence the $0$-faces of $\InjHull(X)$ are precisely
$\bigl(\{f_\lambda\}\bigr)_{\lambda\in\mathbb Y_N}$.
\end{prf}

\begin{example} Here is an example that illustrates the formula
(\ref{numberofsites}) with $N=9$ for the partition $\lambda=(5,3,3,2)$.
Each site in $\mathfrak X$ is represented thrice in the picture
(in other words, it shows three fundamental domains in the universal cover).
The displayed numbers are the values of $\tilde f_\lambda(L)$. (The dashed
rectangle will serve as an illustration of Remark~\ref{inaffhypplane}.)
\setlength{\unitlength}{0.00057in}
\begin{center}
\begin{picture}(10824,3339)(0,-10)
\path(12,3012)(312,3312)(612,3012)
	(912,3312)(1212,3012)(1512,3312)
	(1812,3012)(2112,3312)(2412,3012)
	(2712,3312)(3012,3012)(3312,3312)
	(3612,3012)(3912,3312)(4212,3012)
	(4512,3312)(4812,3012)(5112,3312)
	(5412,3012)(5712,3312)(6012,3012)
	(6312,3312)(6612,3012)(6912,3312)
	(7212,3012)(7512,3312)(7812,3012)
	(8112,3312)(8412,3012)(8712,3312)
	(9012,3012)(9312,3312)(9612,3012)
	(9912,3312)(10212,3012)(10512,3312)
	(10812,3012)(7812,12)(7512,312)
	(7212,12)(6912,312)(6612,12)
	(6312,312)(6012,12)(5712,312)
	(5412,12)(5112,312)(4812,12)
	(4512,312)(4212,12)(3912,312)
	(3612,12)(3312,312)(3012,12)(12,3012)
\texture{55888888 88555555 5522a222 a2555555 55888888 88555555 552a2a2a 2a555555 
	55888888 88555555 55a222a2 22555555 55888888 88555555 552a2a2a 2a555555 
	55888888 88555555 5522a222 a2555555 55888888 88555555 552a2a2a 2a555555 
	55888888 88555555 55a222a2 22555555 55888888 88555555 552a2a2a 2a555555 }
\shade\path(1512,1512)(2412,2412)(2712,2112)
	(3012,2412)(3612,1812)(4212,2412)
	(5112,1512)(5412,1812)(5712,1512)
	(6312,2112)(6912,1512)(7812,2412)
	(8112,2112)(8412,2412)(9012,1812)
	(9612,2412)(9912,2112)(9012,1212)
	(8412,1812)(8112,1512)(7812,1812)
	(6912,912)(6312,1512)(5712,912)
	(5412,1212)(5112,912)(4212,1812)
	(3612,1212)(3012,1812)(2712,1512)
	(2412,1812)(1812,1212)(1512,1512)
\path(1512,1512)(2412,2412)(2712,2112)
	(3012,2412)(3612,1812)(4212,2412)
	(5112,1512)(5412,1812)(5712,1512)
	(6312,2112)(6912,1512)(7812,2412)
	(8112,2112)(8412,2412)(9012,1812)
	(9612,2412)(9912,2112)(9012,1212)
	(8412,1812)(8112,1512)(7812,1812)
	(6912,912)(6312,1512)(5712,912)
	(5412,1212)(5112,912)(4212,1812)
	(3612,1212)(3012,1812)(2712,1512)
	(2412,1812)(1812,1212)(1512,1512)
\thicklines
\path(1812,1212)(3012,12)(4512,1512)(4212,1812)(3612,1212)(3012,1812)(2712,1512)
	(2412,1812)(1812,1212)
\dashline[24]{80}(3912,2112)(5112,3312)(7212,1212)(6012,12)(3912,2112)
\put(1812,1512){\makebox(0,0){\footnotesize$0$}}
\put(4512,1812){\makebox(0,0){\footnotesize$0$}}
\put(7212,1512){\makebox(0,0){\footnotesize$0$}}
\put(2112,1212){\makebox(0,0){\footnotesize$2$}}
\put(2412,912){\makebox(0,0){\footnotesize$5$}}
\put(2712,612){\makebox(0,0){\footnotesize$8$}}
\put(3012,312){\makebox(0,0){\footnotesize$13$}}
\put(2412,1512){\makebox(0,0){\footnotesize$1$}}
\put(2712,1212){\makebox(0,0){\footnotesize$3$}}
\put(3012,912){\makebox(0,0){\footnotesize$5$}}
\put(3312,612){\makebox(0,0){\footnotesize$9$}}
\put(3612,312){\makebox(0,0){\footnotesize$13$}}
\put(3012,1512){\makebox(0,0){\footnotesize$1$}}
\put(3312,1212){\makebox(0,0){\footnotesize$2$}}
\put(3612,912){\makebox(0,0){\footnotesize$5$}}
\put(3912,612){\makebox(0,0){\footnotesize$8$}}
\put(4212,312){\makebox(0,0){\footnotesize$11$}}
\put(3912,1212){\makebox(0,0){\footnotesize$2$}}
\put(4212,912){\makebox(0,0){\footnotesize$4$}}
\put(4512,612){\makebox(0,0){\footnotesize$6$}}
\put(4812,312){\makebox(0,0){\footnotesize$9$}}
\put(4212,1512){\makebox(0,0){\footnotesize$1$}}
\put(4512,1212){\makebox(0,0){\footnotesize$2$}}
\put(4812,912){\makebox(0,0){\footnotesize$3$}}
\put(5112,612){\makebox(0,0){\footnotesize$5$}}
\put(5412,312){\makebox(0,0){\footnotesize$9$}}
\put(5412,912){\makebox(0,0){\footnotesize$1$}}
\put(5712,612){\makebox(0,0){\footnotesize$4$}}
\put(6012,312){\makebox(0,0){\footnotesize$7$}}
\put(6012,912){\makebox(0,0){\footnotesize$2$}}
\put(6312,612){\makebox(0,0){\footnotesize$4$}}
\put(6612,312){\makebox(0,0){\footnotesize$9$}}
\put(6312,1212){\makebox(0,0){\footnotesize$1$}}
\put(6612,912){\makebox(0,0){\footnotesize$2$}}
\put(6912,612){\makebox(0,0){\footnotesize$6$}}
\put(7212,312){\makebox(0,0){\footnotesize$11$}}
\put(7212,912){\makebox(0,0){\footnotesize$3$}}
\put(7512,612){\makebox(0,0){\footnotesize$7$}}
\put(7812,312){\makebox(0,0){\footnotesize$11$}}
\put(7512,1212){\makebox(0,0){\footnotesize$2$}}
\put(7812,912){\makebox(0,0){\footnotesize$5$}}
\put(8112,612){\makebox(0,0){\footnotesize$8$}}
\put(7812,1512){\makebox(0,0){\footnotesize$1$}}
\put(8112,1212){\makebox(0,0){\footnotesize$3$}}
\put(8412,912){\makebox(0,0){\footnotesize$5$}}
\put(8412,1512){\makebox(0,0){\footnotesize$1$}}
\put(8712,1212){\makebox(0,0){\footnotesize$2$}}
\put(4812,2112){\makebox(0,0){\footnotesize$2$}}
\put(5112,2412){\makebox(0,0){\footnotesize$5$}}
\put(5412,2712){\makebox(0,0){\footnotesize$8$}}
\put(5712,3012){\makebox(0,0){\footnotesize$13$}}
\put(6012,2712){\makebox(0,0){\footnotesize$9$}}
\put(6312,2412){\makebox(0,0){\footnotesize$5$}}
\put(6612,2112){\makebox(0,0){\footnotesize$2$}}
\put(6912,1812){\makebox(0,0){\footnotesize$1$}}
\put(6012,2112){\makebox(0,0){\footnotesize$2$}}
\put(5712,2412){\makebox(0,0){\footnotesize$5$}}
\put(5112,3012){\makebox(0,0){\footnotesize$11$}}
\put(4812,2712){\makebox(0,0){\footnotesize$7$}}
\put(4512,2412){\makebox(0,0){\footnotesize$3$}}
\put(5112,1812){\makebox(0,0){\footnotesize$1$}}
\put(5712,1812){\makebox(0,0){\footnotesize$1$}}
\put(5412,2112){\makebox(0,0){\footnotesize$3$}}
\end{picture}
\end{center}
\end{example}

\begin{remark}\label{inaffhypplane}
For $(j,k)\in\mathcal L_\lambda$ (with $0\leqslant j\leqslant k<N$) the formula
\begin{equation}\label{complementdiscrete}
f_\lambda(j)+f_\lambda(k)=d(j,k)
\end{equation}
from Theorem~\ref{zerocells} (see also Definition~\ref{HLdjk})
has the following interpretation: The shapes of
$\alpha_\lambda(L)|_{L=(j,j)=(j,j+N)}$ and $\alpha_\lambda(L)|_{L=(k,k)}$
together with the rim $\mathcal L_\lambda$ emanating from the position
$(j,k)$ and ending at the position $(k,j+N)$
\mbox{$\bigl({}=(j,k)$} as a site in \mbox{$\mathfrak X\bigr)$}
tessellate the $(N-(k-j)+1)$\/$\times$\/$(k-j+1)$ rectangle with corners $(j,k)$,
$(j,j+N)$, $(k,j+N)$, $(k,k)$. Hence
\begin{align*}
f_\lambda(j)+f_\lambda(k)&=\bigl|\alpha_\lambda(L)|_{L=(j,j)=(j,j+N)}\bigr|+
\bigl|\alpha_\lambda(L)|_{L=(k,k)}\bigr|\\
&=(N-(k-j)+1)(k-j+1)-(N+1)=d(j,k).
\end{align*}
As an illustration we look at the example above and take $(j,k)=(5,8)$.
The picture shows the corresponding dashed $7\times4$ rectangle. So
$$f_\lambda(5)+f_\lambda(8)=7+11=7\cdot4-10.$$
The formula (\ref{complement}) in Section~\ref{continuous} is a continuous
analogue of (\ref{complementdiscrete}).
\end{remark}

\begin{theorem}[S.~Herrmann, M.~Joswig \cite{HJ}]\label{HJThm}
The number of $v$-faces in an injective hull of a metric space with $N$ points
is at most
$$2^{N-2v-1}\frac{N}{N-v}\binom{N-v}{v},$$
and for each $N$ there is a metric space attaining those upper bounds
uniformly for all $v$.
\end{theorem}

Let $\lambda\in\mathbb Y_N$ with outer rim $\mathcal L_\lambda$ as in
Proposition~\ref{loopsandpartitions}. Let us enumerate those sites
$L_1=(j_1+1,k_1),\dots,L_s=(j_s+1,k_s)$ in the outer rim
$\mathcal L_\lambda$ that correspond to the inner corners
$(j_1,k_1+1),\dots,(j_s,k_s+1)$ of $\lambda$.
The partition $\lambda$ covers exactly $s$ partitions, namely those partitions
$\mu$ that are got by removing from $\lambda$ one of its inner corners.
The $2^s$ partitions $\nu$ that are got from $\lambda$ by removing $i$ (running
from $0$ to $s$) of its inner corners constitute a Boolean lattice.
We shall recognize that the convex hull of those $2^s$ vertices
$f_\nu$ belongs to the injective hull $\InjHull(X)$.
In fact, let $f=\sum_\nu a_\nu f_\nu$ be such a convex combination. Then from
$\tilde f_\nu(L)\geqslant0$ for all $L\in\mathfrak X$ we have
$\tilde f(L)\geqslant0$ for all $L\in\mathfrak X$, that is, $f\in\Delta(X)$.
For each $L\in\mathcal L_\lambda-\bigl\{L_1,\dots,L_s\bigr\}$ we have
$\tilde f_\nu(L)=0$ and hence $\tilde f(L)=0$.
Since $(j_i,k_i),(j_i+1,k_i+1)\in\mathcal L_\lambda-\bigl\{L_1,\dots,L_s\bigr\}$,
we still get (as a generalization of Remark~\ref{exhaustX})
\begin{equation}\label{exhaustXgen}
X=\bigcup_{(j,k)\in\mathcal L_\lambda-\{L_1,\dots,L_s\}}\{j,k\}
\end{equation}
from which we conclude that $f\in\InjHull(X)$. (Note also that
from (\ref{exhaustXgen}) we get $s\leqslant\bigl\lfloor\frac N2\bigr\rfloor$.
This is Develin's bound for the maximal possible dimension of the
injective hull of a metric space with $N$ points (see \cite{De}), a result that
is superseded by Theorem~\ref{HJThm}.)

A systematic count of $v$-faces can be done as follows: Each choice of $v$
among the $s$ inner corners of $\lambda$ results in a $v$-face with
$f_\lambda$ as its ``top vertex''. In this way, we get $\binom sv$ $v$-faces
with ``top vertex'' $f_\lambda$.
The total number of $v$-faces that we get from this procedure
(and we shall soon see that every $v$-face occurs in such a way) is thus
\begin{equation}\label{vfacesinE}
\mbox{\#($v$-faces in $\InjHull(X)$)}=
\sum_{\lambda\in\mathbb Y_N}\binom{s(\lambda)}{v}
\end{equation}
where $s(\lambda)$ denotes the number of inner corners of the
partition $\lambda$.

\begin{lemma}
The number of partitions with first part at most $j$ and having at most $k$
parts that cover exactly $s$ partitions is $a(j,k,s):=
\binom{j}{s}\binom{k}{s}$. 
\end{lemma}
\begin{prf}
We have to count the partitions with first part at most $j$ and at most $k$
parts with exactly $s$ inner corners. This amounts to choosing $s$ integers
$1\leqslant j_1<\dots<j_s\leqslant j$ and $s$ integers
$1\leqslant k_1<\dots<k_s\leqslant k$ because
the pairs $(j_i,k_{s+1-i})$ for $1\leqslant i\leqslant s$ determine the positions
of the $s$ inner corners of a partition in the given range in an evident way
as illustrated in the picture.
\setlength{\unitlength}{0.00048in}
\begin{center}
\begin{picture}(6324,5739)(0,-10)
\put(2112,2712){\blacken\ellipse{120}{120}}
\put(2112,2712){\ellipse{120}{120}}
\put(2712,2712){\blacken\ellipse{120}{120}}
\put(2712,2712){\ellipse{120}{120}}
\put(4212,3012){\blacken\ellipse{120}{120}}
\put(4212,3012){\ellipse{120}{120}}
\put(5112,3312){\blacken\ellipse{120}{120}}
\put(5112,3312){\ellipse{120}{120}}
\path(1812,2712)(2112,2412)(2412,2712)
\path(2412,2712)(2712,2412)(3012,2712)
\path(3912,3012)(4212,2712)(4512,3012)
\path(4812,3312)(5112,3012)
\path(1212,2112)(2112,3012)(2412,2712)
	(2712,3012)(3312,2412)(4212,3312)
	(4512,3012)(5112,3612)(5412,3312)
	(2712,612)(1212,2112)
\path(1212,2112)(612,2712)
\path(1212,2112)(612,2712)
\path(5412,3312)(5712,3612)
\path(5412,3312)(5712,3612)
\path(118.066,2048.360)(12.000,2112.000)(75.640,2005.934)
\path(12,2112)(2112,12)
\path(12,2112)(2112,12)
\path(2005.934,75.640)(2112.000,12.000)(2048.360,118.066)
\path(3375.640,118.066)(3312.000,12.000)(3418.066,75.640)
\path(3312,12)(6312,3012)
\path(3312,12)(6312,3012)
\path(6248.360,2905.934)(6312.000,3012.000)(6205.934,2948.360)
\path(612,2712)(3612,5712)(5712,3612)
\path(1512,1812)(2112,2412)(3312,1212)
\path(1812,1512)(2712,2412)(3612,1512)
\path(2112,1212)(3312,2412)(3912,1812)
\path(2412,912)(4212,2712)(4512,2412)
\path(3612,2712)(4212,2112)
\path(4512,3012)(4812,2712)
\path(1512,2412)(3012,912)
\put(3612,1212){\makebox(0,0){\footnotesize$j_1$}}
\put(3912,1512){\makebox(0,0){\footnotesize$j_2$}}
\put(4812,2412){\makebox(0,0){\footnotesize$\iddots$}}
\put(5412,3012){\makebox(0,0){\footnotesize$j_s$}}
\put(2412,612){\makebox(0,0){\footnotesize$k_1$}}
\put(2112,912){\makebox(0,0){\footnotesize$k_2$}}
\put(1512,1512){\makebox(0,0){\footnotesize$\ddots$}}
\put(1212,1812){\makebox(0,0){\footnotesize$k_s$}}
\put(5112,1212){\makebox(0,0){$j$}}
\put(612,612){\makebox(0,0){$k$}}
\end{picture}
\end{center}
The correspondence is clearly bijective.
\end{prf}
 
\begin{proposition}\label{coverinY}
The number of partitions in $\mathbb Y_N$ that cover exactly $s$ partitions
(that is, with exactly $s$ inner corners) is $\binom{N}{2s}$.
\end{proposition}
\begin{prf}
Let us compute the number in question by inclusion-exclusion and using
the previous lemma and its notation $a(j,k,s)=\binom{j}{s}\binom{k}{s}$
for the number of partitions with first part at most $j$
and having at most $k$ parts that cover exactly $s$ partitions.

The formula is obviously true for $s=0$ (corresponding to the empty partition,
which is contained in every $\mathbb Y_N$); hence we assume that
$s\geqslant1$. The following figure explains the inductive step for the
inclusion-exclusion argument.
\begin{center}
\begin{picture}(8124,5739)(0,-10)
\path(118.066,2048.360)(12.000,2112.000)(75.640,2005.934)
\path(12,2112)(2112,12)
\path(12,2112)(2112,12)
\path(2005.934,75.640)(2112.000,12.000)(2048.360,118.066)
\path(3375.640,118.066)(3312.000,12.000)(3418.066,75.640)
\path(3312,12)(6312,3012)
\path(3312,12)(6312,3012)
\path(6248.360,2905.934)(6312.000,3012.000)(6205.934,2948.360)
\path(2712,612)(7512,5412)(7212,5712)
	(6912,5412)(6612,5712)(6312,5412)
	(6012,5712)(5712,5412)(5412,5712)
	(5112,5412)(4812,5712)(4512,5412)
	(4212,5712)(912,2412)(2712,612)
\path(2712,612)(5712,3612)(3612,5712)
	(612,2712)(2712,612)
\path(6375.640,3118.066)(6312.000,3012.000)(6418.066,3075.640)
\path(6312,3012)(8112,4812)
\path(8048.360,4705.934)(8112.000,4812.000)(8005.934,4748.360)
\put(5112,1212){\makebox(0,0)[l]{$N-k$}}
\put(612,612){\makebox(0,0){$k$}}
\put(7512,3612){\makebox(0,0)[l]{$k-1$}}
\end{picture}
\end{center}
The number $A(N,k,s)$ of partitions in $\mathbb Y_N$
having at most $k$ parts that cover exactly $s$ partitions is
$$A(N,k,s)=A(N,k-1,s)+a(N-k,k,s)-a(N-k,k-1,s).$$
Furthermore, $A(N,1,s)=a(N-1,1,s)$. 
The number $A(N,N-1,s)$ that we want to compute is therefore
\begin{align*}
A(N,N-1,s)&=\sum_{k=1}^{N-1}a(N-k,k,s)-\sum_{k=2}^{N-1}a(N-k,k-1,s)\\
&=\sum_{k=1}^{N-1}\binom{N-k}{s}\binom{k}{s}-\sum_{k=2}^{N-1}\binom{N-k}{s}
\binom{k-1}{s}\\
\intertext{and combining the binomial coefficients $\binom ks-\binom{k-1}{s}
=\binom{k-1}{s-1}$ after adding the zero summand
$-\binom{N-k}{s}\binom{k-1}{s}\bigr|_{k=1}$
(recall that $s\geqslant1$) in the second sum, we get}
&=\sum_{k=1}^{N-1}\binom{N-k}{s}\binom{k-1}{s-1}=\binom{N}{2s},
\end{align*}
where the last equality is clear from the following combinatorial
interpretation: To choose $2s$ integers $1\leqslant i_1<\dots<i_{s-1}
<i_s<i_{s+1}<\dots<i_{2s}\leqslant N$ is tantamount to first fixing $i_s=k\in
\{1,\dots,N-1\}$ and then choosing $s-1$ integers in $\{1,\dots,k-1\}$
and choosing $s$ integers in $\{k+1,\dots,N\}$.
\end{prf} 

\begin{lemma}\label{numberfaces}
The number of $v$-faces in $\InjHull(X_N)$ is (at least and in fact exactly)
$$\mbox{\textup{\#($v$-faces in $\InjHull(X_N)$)}}
=\sum_{s=v}^{\lfloor N/2\rfloor}\binom{N}{2s}\binom sv.$$
\end{lemma}
\begin{prf}
This is clear from Proposition~\ref{coverinY} together with (\ref{vfacesinE}).
\end{prf}

\begin{theorem}\label{allfaces}
Consider the $N$-point metric space $X_N=\{0,1,\dots,N-1\}$ with
metric $d(j,k)=|k-j|(N-|k-j|)$.
The number of $v$-faces in its injective hull $\InjHull(X_N)$ (as in
Definition~\ref{InjHullX}) is
\begin{equation}\label{numbervfaces}
\mbox{\textup{\#($v$-faces in $\InjHull(X_N)$)}}
=2^{N-2v-1}\frac{N}{N-v}\binom{N-v}{v}.
\end{equation}
\end{theorem}
\begin{prf}
For $v=0$ we have already proved this in Theorem~\ref{zerocells} because
$|\mathbb Y_N|=2^{N-1}$. Hence we assume that $v\geqslant1$. Using
Lemma~\ref{numberfaces} and three instances of the binomial series expansion
(as a formal power series)
\begin{equation}\tag{$\ast$}
\sum_{n=a}^\infty\binom na z^n=\frac{z^a}{(1-z)^{a+1}}\qquad%
\mbox{for $a\in\mathbb Z_{\geqslant0}$}
\end{equation}
(namely, twice in the direction from left to right and once in the opposite
direction), we\newpage
\noindent compute the generating series
\begin{align*}
\lefteqn{\sum_{N\geqslant2}\mbox{\#($v$-faces in $\InjHull(X_N)$)}\,q^N}
\qquad\\[-2mm]
&\stackrel{\phantom{(\ast)}}{=}\sum_{N\geqslant2}
\sum_{s=v}^{\lfloor N/2\rfloor}\binom{N}{2s}\binom sv q^N
=\sum_{s=v}^\infty\sum_{N=2s}^\infty\binom{N}{2s}\binom sv q^N
=\sum_{s=v}^\infty\binom sv\sum_{N=2s}^\infty\binom{N}{2s}q^N\\
&\stackrel{(\ast)}{=}
\sum_{s=v}^\infty\binom sv\frac{q^{2s}}{(1-q)^{2s+1}}=\frac{1}{(1-q)}
\sum_{s=v}^\infty\binom sv\left(\frac{q^2}{(1-q)^2}\right)^s\\
&\stackrel{(\ast)}{=}
\frac{1}{(1-q)}\frac{\left(\dfrac{q^2}{(1-q)^2}\right)^v}{
\left(1-\dfrac{q^2}{(1-q)^2}\right)^{v+1}}
=\frac{(1-q)q^{2v}}{(1-2q)^{v+1}}=(1-q)q^v 2^{-v}\frac{(2q)^v}{(1-2q)^{v+1}}\\
&\stackrel{(\ast)}{=}(1-q)q^v 2^{-v}\sum_{n=v}^\infty\binom nv (2q)^n\\
&\stackrel{\phantom{(\ast)}}{=}q^v 2^{-v-1}\sum_{N=2v}^\infty2\binom{N-v}{v}(2q)^{N-v}
-q^{v+1} 2^{-v}\sum_{N=2v+1}^\infty\binom{N-v-1}{v}(2q)^{N-v-1}\\
\intertext{and using $2\binom{N-v}{v}-\binom{N-v-1}{v}=
\frac{N}{N-v}\binom{N-v}{v}$ and $\binom{N-v-1}{v}\bigr|_{N=2v}=0$, we get
finally}
&\stackrel{\phantom{(\ast)}}{=}
\sum_{N=2v}^\infty 2^{N-2v-1}\frac{N}{N-v}\binom{N-v}{v}\,q^N.
\end{align*}
Hence the formula (\ref{numbervfaces}) is proved. According to
Theorem~\ref{HJThm}, $X_N$ is a metric space such that $\InjHull(X_N)$ attains
the maximal possible number of $v$-faces for all $v$.
\end{prf}

\begin{remark}\label{ProofSummaryThm}
Theorem~\ref{SummaryThm} from the introductory section follows as a
corollary. In fact, the vertices of $\InjHull(X_N)$ were described in
Theorem~\ref{zerocells} as $f_\lambda\in\InjHull(X_N)$ for $\lambda\in
\mathbb Y_N$. The considerations after Theorem~\ref{HJThm} showed that
if $\mathbb Y_N\ni\lambda$ covers $\mu$, then the segment with endpoints
$f_\lambda$ and $f_\mu$ is an edge in $\InjHull(X_N)$. In this way we get
$2^{N-3}N$ edges in $\InjHull(X_N)$, and according to
Theorem~\ref{allfaces}, those $2^{N-3}N$ edges are all the $1$-faces in
$\InjHull(X_N)$.
\end{remark}

\enlargethispage*{5mm}
\begin{remark}
For $v=1$ we get $2^{N-3}N$ edges in $\InjHull(X_N)$ and thus in
$\operatorname{Hasse}(\mathbb Y_N)$. This result generalizes to the
fact that the Hasse diagram of the poset of abelian ideals in a Borel
subalgebra of a complex simple Lie algebra of rank $n$ has $2^{n-2}(n+1)$
edges (see \cite[Theorem~4.1]{Pa}).

For $v\geqslant2$ such a simple (only rank-dependent) census breaks down.
Recall that $\InjHull(X_5)$ (type $\mathsf A_4$)
has five $2$-faces and no higher-dimensional faces.
But the poset of abelian ideals in a Borel subalgebra of
$\mathfrak{so}(8,\mathbb C)$ (type $\mathsf D_4$) contains a Boolean subposet
of rank $3$: If $\alpha_1,\alpha_2,\alpha_3,\alpha_4$
are the simple roots and $\theta=\alpha_1+2\alpha_2+\alpha_3+\alpha_4$ is
the highest root, then the eight ideals in question are
$$\mathfrak g_\theta\oplus\mathfrak g_{\theta-\alpha_2}
\oplus\left(\begin{array}{@{}c@{}}
\mathfrak g_{\theta-\alpha_2-\alpha_1}\\\textup{or}\\0\end{array}\right)
\oplus\left(\begin{array}{@{}c@{}}
\mathfrak g_{\theta-\alpha_2-\alpha_3}\\\textup{or}\\0\end{array}\right)
\oplus\left(\begin{array}{@{}c@{}}
\mathfrak g_{\theta-\alpha_2-\alpha_4}\\\textup{or}\\0\end{array}\right),$$
where $\mathfrak g_\varphi$ denotes the root space for the
root $\varphi$.
Hence the Hasse diagram contains the $1$-skeleton of a cube,
which has one $3$-face and six $2$-faces.
\end{remark}

\newpage\enlargethispage*{1cm}
\begin{example}[The Hasse diagram for $N=5$] Its vertices
$\lambda\in\mathbb Y_5$ are encoded by the graphs with vertex set $X_5$, and
there is an edge between $j$ and $k$ if and only if
$(j,k)\in\mathcal L_\lambda$ (the outer rim of $\lambda$ as defined in
Proposition~\ref{loopsandpartitions}).
%0 (3383,2042)
%1 (669,1160)
%2 (2347,3469)
%3 (2347,615)
%4 (669,2924)
\setlength{\unitlength}{0.00017in}
\begin{center}
\xymatrix@R=12mm@C=1mm{&&&{\begin{picture}(3991,4089)(0,-10)
\thicklines
\put(2439,3754){\ellipse{600}{600}} %2
\put(3683,2042){\makebox(0,0){\tiny$0$}}
\put(427,984){\makebox(0,0){\tiny$1$}}
\put(2439,3754){\makebox(0,0){\tiny$2$}}
\put(2439,330){\makebox(0,0){\tiny$3$}}
\put(427,3100){\makebox(0,0){\tiny$4$}}
\path(2347,3469)(3383,2042)
\path(2347,3469)(669,1160)
\path(2347,3469)(2347,615)
\path(2347,3469)(669,2924)
\end{picture}}&&
{\begin{picture}(3991,4089)(0,-10)
\thicklines
\put(2439,330){\ellipse{600}{600}} %3
\put(3683,2042){\makebox(0,0){\tiny$0$}}
\put(427,984){\makebox(0,0){\tiny$1$}}
\put(2439,3754){\makebox(0,0){\tiny$2$}}
\put(2439,330){\makebox(0,0){\tiny$3$}}
\put(427,3100){\makebox(0,0){\tiny$4$}}
\path(2347,615)(3383,2042)
\path(2347,615)(669,1160)
\path(2347,615)(2347,3469)
\path(2347,615)(669,2924)
\end{picture}}\\
&&&{\begin{picture}(3991,4089)(0,-10)
\thicklines
\put(3683,2042){\makebox(0,0){\tiny$0$}}
\put(427,984){\makebox(0,0){\tiny$1$}}
\put(2439,3754){\makebox(0,0){\tiny$2$}}
\put(2439,330){\makebox(0,0){\tiny$3$}}
\put(427,3100){\makebox(0,0){\tiny$4$}}
\path(2347,3469)(3383,2042)
\path(2347,3469)(669,1160)
\path(2347,3469)(2347,615)
\path(2347,3469)(669,2924)
\path(669,1160)(2347,615)
\end{picture}}\ar@{-}[u]&&{\begin{picture}(3991,4089)(0,-10)
\thicklines
\put(3683,2042){\makebox(0,0){\tiny$0$}}
\put(427,984){\makebox(0,0){\tiny$1$}}
\put(2439,3754){\makebox(0,0){\tiny$2$}}
\put(2439,330){\makebox(0,0){\tiny$3$}}
\put(427,3100){\makebox(0,0){\tiny$4$}}
\path(2347,615)(3383,2042)
\path(2347,615)(669,1160)
\path(2347,615)(2347,3469)
\path(2347,615)(669,2924)
\path(2347,3469)(669,2924)
\end{picture}}\ar@{-}[u]\\
{\begin{picture}(3991,4089)(0,-10)
\thicklines
\put(427,984){\ellipse{600}{600}} %1
\put(3683,2042){\makebox(0,0){\tiny$0$}}
\put(427,984){\makebox(0,0){\tiny$1$}}
\put(2439,3754){\makebox(0,0){\tiny$2$}}
\put(2439,330){\makebox(0,0){\tiny$3$}}
\put(427,3100){\makebox(0,0){\tiny$4$}}
\path(669,1160)(3383,2042)
\path(669,1160)(2347,3469)
\path(669,1160)(2347,615)
\path(669,1160)(669,2924)
\end{picture}}&&{\begin{picture}(3991,4089)(0,-10)
\thicklines
\put(3683,2042){\makebox(0,0){\tiny$0$}}
\put(427,984){\makebox(0,0){\tiny$1$}}
\put(2439,3754){\makebox(0,0){\tiny$2$}}
\put(2439,330){\makebox(0,0){\tiny$3$}}
\put(427,3100){\makebox(0,0){\tiny$4$}}
\path(3383,2042)(2347,3469)
\path(669,1160)(2347,3469)
\path(669,1160)(2347,615)
\path(669,1160)(669,2924)
\path(2347,3469)(669,2924)
\end{picture}}\ar@{-}[ur]&
{\raisebox{-1.1em}[0pt][0pt]{\fbox{\begin{picture}(3991,4089)(0,-10)
\thicklines
\put(3683,2042){\makebox(0,0){\tiny$0$}}
\put(427,984){\makebox(0,0){\tiny$1$}}
\put(2439,3754){\makebox(0,0){\tiny$2$}}
\put(2439,330){\makebox(0,0){\tiny$3$}}
\put(427,3100){\makebox(0,0){\tiny$4$}}
\path(669,1160)(2347,615)
\path(3383,2042)(2347,3469)
\path(2347,3469)(669,2924)
\end{picture}}}}&{\begin{picture}(3991,4089)(0,-10)
\thicklines
\put(3683,2042){\makebox(0,0){\tiny$0$}}
\put(427,984){\makebox(0,0){\tiny$1$}}
\put(2439,3754){\makebox(0,0){\tiny$2$}}
\put(2439,330){\makebox(0,0){\tiny$3$}}
\put(427,3100){\makebox(0,0){\tiny$4$}}
\path(3383,2042)(2347,615)
\path(669,1160)(2347,615)
\path(2347,3469)(2347,615)
\path(2347,3469)(669,2924)
\path(3383,2042)(2347,3469)
\end{picture}}\ar@{-}[ul]\ar@{-}[ur]
&{\raisebox{-1.1em}[0pt][0pt]{\fbox{\begin{picture}(3991,4089)(0,-10)
\thicklines
\put(3683,2042){\makebox(0,0){\tiny$0$}}
\put(427,984){\makebox(0,0){\tiny$1$}}
\put(2439,3754){\makebox(0,0){\tiny$2$}}
\put(2439,330){\makebox(0,0){\tiny$3$}}
\put(427,3100){\makebox(0,0){\tiny$4$}}
\path(2347,3469)(669,2924)
\path(3383,2042)(2347,615)
\path(669,1160)(2347,615)
\end{picture}}}}&{\begin{picture}(3991,4089)(0,-10)
\thicklines
\put(3683,2042){\makebox(0,0){\tiny$0$}}
\put(427,984){\makebox(0,0){\tiny$1$}}
\put(2439,3754){\makebox(0,0){\tiny$2$}}
\put(2439,330){\makebox(0,0){\tiny$3$}}
\put(427,3100){\makebox(0,0){\tiny$4$}}
\path(3383,2042)(2347,615)
\path(669,1160)(2347,615)
\path(669,1160)(669,2924)
\path(2347,3469)(669,2924)
\path(2347,615)(669,2924)
\end{picture}}\ar@{-}[ul]&&{\begin{picture}(3991,4089)(0,-10)
\thicklines
\put(427,3100){\ellipse{600}{600}} %4
\put(3683,2042){\makebox(0,0){\tiny$0$}}
\put(427,984){\makebox(0,0){\tiny$1$}}
\put(2439,3754){\makebox(0,0){\tiny$2$}}
\put(2439,330){\makebox(0,0){\tiny$3$}}
\put(427,3100){\makebox(0,0){\tiny$4$}}
\path(669,2924)(3383,2042)
\path(669,2924)(669,1160)
\path(669,2924)(2347,3469)
\path(669,2924)(2347,615)
\end{picture}}\\
&&{\begin{picture}(3991,4089)(0,-10)
\thicklines
\put(3683,2042){\makebox(0,0){\tiny$0$}}
\put(427,984){\makebox(0,0){\tiny$1$}}
\put(2439,3754){\makebox(0,0){\tiny$2$}}
\put(2439,330){\makebox(0,0){\tiny$3$}}
\put(427,3100){\makebox(0,0){\tiny$4$}}
\path(3383,2042)(2347,3469)
\path(669,1160)(2347,3469)
\path(669,1160)(669,2924)
\path(669,1160)(2347,615)
\path(3383,2042)(669,1160)
\end{picture}}\ar@{-}[ull]\ar@{-}[u]&
{\raisebox{-2.6em}[0pt][0pt]{\fbox{\begin{picture}(3991,4089)(0,-10)
\thicklines
\put(3683,2042){\makebox(0,0){\tiny$0$}}
\put(427,984){\makebox(0,0){\tiny$1$}}
\put(2439,3754){\makebox(0,0){\tiny$2$}}
\put(2439,330){\makebox(0,0){\tiny$3$}}
\put(427,3100){\makebox(0,0){\tiny$4$}}
\path(3383,2042)(2347,3469)
\path(669,1160)(2347,615)
\path(669,1160)(669,2924)
\end{picture}}}}&{\begin{picture}(3991,4089)(0,-10)
\thicklines
\put(3683,2042){\makebox(0,0){\tiny$0$}}
\put(427,984){\makebox(0,0){\tiny$1$}}
\put(2439,3754){\makebox(0,0){\tiny$2$}}
\put(2439,330){\makebox(0,0){\tiny$3$}}
\put(427,3100){\makebox(0,0){\tiny$4$}}
\path(3383,2042)(2347,615)
\path(669,1160)(2347,615)
\path(669,1160)(669,2924)
\path(2347,3469)(669,2924)
\path(3383,2042)(2347,3469)
\end{picture}}\ar@{-}[ull]
\ar@{-}[u]\ar@{-}[urr]&
{\raisebox{-2.6em}[0pt][0pt]{\fbox{\begin{picture}(3991,4089)(0,-10)
\thicklines
\put(3683,2042){\makebox(0,0){\tiny$0$}}
\put(427,984){\makebox(0,0){\tiny$1$}}
\put(2439,3754){\makebox(0,0){\tiny$2$}}
\put(2439,330){\makebox(0,0){\tiny$3$}}
\put(427,3100){\makebox(0,0){\tiny$4$}}
\path(3383,2042)(2347,615)
\path(669,1160)(669,2924)
\path(2347,3469)(669,2924)
\end{picture}}}}&{\begin{picture}(3991,4089)(0,-10)
\thicklines
\put(3683,2042){\makebox(0,0){\tiny$0$}}
\put(427,984){\makebox(0,0){\tiny$1$}}
\put(2439,3754){\makebox(0,0){\tiny$2$}}
\put(2439,330){\makebox(0,0){\tiny$3$}}
\put(427,3100){\makebox(0,0){\tiny$4$}}
\path(3383,2042)(669,2924)
\path(669,1160)(669,2924)
\path(2347,3469)(669,2924)
\path(2347,615)(669,2924)
\path(3383,2042)(2347,615)
\end{picture}}\ar@{-}[u]\ar@{-}[urr]\\
&&&{\begin{picture}(3991,4089)(0,-10)
\thicklines
\put(3683,2042){\makebox(0,0){\tiny$0$}}
\put(427,984){\makebox(0,0){\tiny$1$}}
\put(2439,3754){\makebox(0,0){\tiny$2$}}
\put(2439,330){\makebox(0,0){\tiny$3$}}
\put(427,3100){\makebox(0,0){\tiny$4$}}
\path(3383,2042)(2347,615)
\path(669,1160)(2347,615)
\path(669,1160)(669,2924)
\path(3383,2042)(669,1160)
\path(3383,2042)(2347,3469)
\end{picture}}\ar@{-}[ul]\ar@{-}[ur]
&{\raisebox{-1.9em}[0pt][0pt]{\fbox{\begin{picture}(3991,4089)(0,-10)
\thicklines
\put(3683,2042){\makebox(0,0){\tiny$0$}}
\put(427,984){\makebox(0,0){\tiny$1$}}
\put(2439,3754){\makebox(0,0){\tiny$2$}}
\put(2439,330){\makebox(0,0){\tiny$3$}}
\put(427,3100){\makebox(0,0){\tiny$4$}}
\path(3383,2042)(2347,3469)
\path(3383,2042)(2347,615)
\path(669,1160)(669,2924)
\end{picture}}}}&{\begin{picture}(3991,4089)(0,-10)
\thicklines
\put(3683,2042){\makebox(0,0){\tiny$0$}}
\put(427,984){\makebox(0,0){\tiny$1$}}
\put(2439,3754){\makebox(0,0){\tiny$2$}}
\put(2439,330){\makebox(0,0){\tiny$3$}}
\put(427,3100){\makebox(0,0){\tiny$4$}}
\path(3383,2042)(669,2924)
\path(669,1160)(669,2924)
\path(2347,3469)(669,2924)
\path(3383,2042)(2347,3469)
\path(3383,2042)(2347,615)
\end{picture}}\ar@{-}[ul]\ar@{-}[ur]\\
&&&&{\begin{picture}(3991,4089)(0,-10)
\thicklines
\put(3683,2042){\makebox(0,0){\tiny$0$}}
\put(427,984){\makebox(0,0){\tiny$1$}}
\put(2439,3754){\makebox(0,0){\tiny$2$}}
\put(2439,330){\makebox(0,0){\tiny$3$}}
\put(427,3100){\makebox(0,0){\tiny$4$}}
\path(3383,2042)(669,2924)
\path(669,1160)(669,2924)
\path(3383,2042)(669,1160)
\path(3383,2042)(2347,3469)
\path(3383,2042)(2347,615)
\end{picture}}\ar@{-}[ul]\ar@{-}[ur]\\
&&&&\begin{picture}(3991,4089)(0,-10)
\thicklines
\put(3683,2042){\ellipse{600}{600}} %0
\put(3683,2042){\makebox(0,0){\tiny$0$}}
\put(427,984){\makebox(0,0){\tiny$1$}}
\put(2439,3754){\makebox(0,0){\tiny$2$}}
\put(2439,330){\makebox(0,0){\tiny$3$}}
\put(427,3100){\makebox(0,0){\tiny$4$}}
\path(3383,2042)(669,1160)
\path(3383,2042)(2347,3469)
\path(3383,2042)(2347,615)
\path(3383,2042)(669,2924)
\end{picture}{}\ar@{-}[u]}
\end{center}
The five graphs in the framed boxes encode the $2$-faces of $\InjHull(X_5)$ in
an evident manner (and the twenty edges of $\mathbb Y_5$ or $1$-faces
of $\InjHull(X_5)$ are encoded by the obvious (but not separately displayed)
spanning trees).
\end{example}

\newpage\enlargethispage*{1cm}
\begin{example}[The central cube for $N=6$]
The vertices in the first visualization show the outer rims
(as defined in Proposition~\ref{loopsandpartitions}).
\setlength{\unitlength}{0.0002in}
\begin{center}
\mbox{\xymatrix@R=16mm{&
{\mbox{\begin{picture}(4224,2139)(0,-10)
\path(12,2112)(2112,12)(4212,2112)
\texture{55888888 88555555 5522a222 a2555555 55888888 88555555 552a2a2a 2a555555 
	55888888 88555555 55a222a2 22555555 55888888 88555555 552a2a2a 2a555555 
	55888888 88555555 5522a222 a2555555 55888888 88555555 552a2a2a 2a555555 
	55888888 88555555 55a222a2 22555555 55888888 88555555 552a2a2a 2a555555 }
\shade\path(912,1212)(1212,912)(1512,1212)
	(1812,912)(2112,1212)(2412,912)
	(2712,1212)(3012,912)(3312,1212)
	(2712,1812)(2412,1512)(2112,1812)
	(1812,1512)(1512,1812)(912,1212)
\path(912,1212)(1212,912)(1512,1212)
	(1812,912)(2112,1212)(2412,912)
	(2712,1212)(3012,912)(3312,1212)
	(2712,1812)(2412,1512)(2112,1812)
	(1812,1512)(1512,1812)(912,1212)
\end{picture}}}\ar@{-}[dl]\ar@{-}[d]\ar@{-}[dr]\\
{\mbox{\begin{picture}(4224,2139)(0,-10)
\path(12,2112)(2112,12)(4212,2112)
\path(612,1512)(912,1812)(1212,1512)
	(912,1212)(612,1512)
\texture{55888888 88555555 5522a222 a2555555 55888888 88555555 552a2a2a 2a555555 
	55888888 88555555 55a222a2 22555555 55888888 88555555 552a2a2a 2a555555 
	55888888 88555555 5522a222 a2555555 55888888 88555555 552a2a2a 2a555555 
	55888888 88555555 55a222a2 22555555 55888888 88555555 552a2a2a 2a555555 }
\shade\path(912,1212)(1212,912)(1512,1212)
	(1812,912)(2112,1212)(2712,612)
	(3312,1212)(3012,1512)(2712,1212)
	(2112,1812)(1812,1512)(1512,1812)(912,1212)
\path(912,1212)(1212,912)(1512,1212)
	(1812,912)(2112,1212)(2712,612)
	(3312,1212)(3012,1512)(2712,1212)
	(2112,1812)(1812,1512)(1512,1812)(912,1212)
\end{picture}}}\ar@{-}[d]\ar@{-}[dr]&
{\mbox{\begin{picture}(4224,2139)(0,-10)
\path(12,2112)(2112,12)(4212,2112)
\texture{55888888 88555555 5522a222 a2555555 55888888 88555555 552a2a2a 2a555555 
	55888888 88555555 55a222a2 22555555 55888888 88555555 552a2a2a 2a555555 
	55888888 88555555 5522a222 a2555555 55888888 88555555 552a2a2a 2a555555 
	55888888 88555555 55a222a2 22555555 55888888 88555555 552a2a2a 2a555555 }
\shade\path(912,1212)(1212,912)(1512,1212)
	(2112,612)(2712,1212)(3012,912)
	(3312,1212)(2712,1812)(2112,1212)
	(1512,1812)(912,1212)
\path(912,1212)(1212,912)(1512,1212)
	(2112,612)(2712,1212)(3012,912)
	(3312,1212)(2712,1812)(2112,1212)
	(1512,1812)(912,1212)
\end{picture}}}\ar@{-}[dl]\ar@{-}[dr]&
{\mbox{\begin{picture}(4224,2139)(0,-10)
\path(12,2112)(2112,12)(4212,2112)
\path(3612,1512)(3312,1812)(3012,1512)
	(3312,1212)(3612,1512)
\texture{55888888 88555555 5522a222 a2555555 55888888 88555555 552a2a2a 2a555555 
	55888888 88555555 55a222a2 22555555 55888888 88555555 552a2a2a 2a555555 
	55888888 88555555 5522a222 a2555555 55888888 88555555 552a2a2a 2a555555 
	55888888 88555555 55a222a2 22555555 55888888 88555555 552a2a2a 2a555555 }
\shade\path(912,1212)(1512,612)(2112,1212)
	(2412,912)(2712,1212)(3012,912)
	(3312,1212)(2712,1812)(2412,1512)
	(2112,1812)(1512,1212)(1212,1512)(912,1212)
\path(912,1212)(1512,612)(2112,1212)
	(2412,912)(2712,1212)(3012,912)
	(3312,1212)(2712,1812)(2412,1512)
	(2112,1812)(1512,1212)(1212,1512)(912,1212)
\end{picture}}}\ar@{-}[dl]\ar@{-}[d]\\
{\mbox{\begin{picture}(4224,2139)(0,-10)
\path(12,2112)(2112,12)(4212,2112)
\path(612,1512)(912,1812)(1212,1512)
	(912,1212)(612,1512)
\texture{55888888 88555555 5522a222 a2555555 55888888 88555555 552a2a2a 2a555555 
	55888888 88555555 55a222a2 22555555 55888888 88555555 552a2a2a 2a555555 
	55888888 88555555 5522a222 a2555555 55888888 88555555 552a2a2a 2a555555 
	55888888 88555555 55a222a2 22555555 55888888 88555555 552a2a2a 2a555555 }
\shade\path(912,1212)(1212,912)(1512,1212)
	(2112,612)(2412,912)(2712,612)
	(3312,1212)(3012,1512)(2712,1212)
	(2412,1512)(2112,1212)(1512,1812)(912,1212)
\path(912,1212)(1212,912)(1512,1212)
	(2112,612)(2412,912)(2712,612)
	(3312,1212)(3012,1512)(2712,1212)
	(2412,1512)(2112,1212)(1512,1812)(912,1212)
\end{picture}}}\ar@{-}[dr]&
{\mbox{\begin{picture}(4224,2139)(0,-10)
\path(12,2112)(2112,12)(4212,2112)
\path(612,1512)(912,1812)(1212,1512)
	(912,1212)(612,1512)
\path(3612,1512)(3312,1812)(3012,1512)
	(3312,1212)(3612,1512)
\texture{55888888 88555555 5522a222 a2555555 55888888 88555555 552a2a2a 2a555555 
	55888888 88555555 55a222a2 22555555 55888888 88555555 552a2a2a 2a555555 
	55888888 88555555 5522a222 a2555555 55888888 88555555 552a2a2a 2a555555 
	55888888 88555555 55a222a2 22555555 55888888 88555555 552a2a2a 2a555555 }
\shade\path(912,1212)(1512,612)(2112,1212)
	(2712,612)(3312,1212)(3012,1512)
	(2712,1212)(2112,1812)(1512,1212)
	(1212,1512)(912,1212)
\path(912,1212)(1512,612)(2112,1212)
	(2712,612)(3312,1212)(3012,1512)
	(2712,1212)(2112,1812)(1512,1212)
	(1212,1512)(912,1212)
\end{picture}}}\ar@{-}[d]&
{\mbox{\begin{picture}(4224,2139)(0,-10)
\path(12,2112)(2112,12)(4212,2112)
\path(3612,1512)(3312,1812)(3012,1512)
	(3312,1212)(3612,1512)
\texture{55888888 88555555 5522a222 a2555555 55888888 88555555 552a2a2a 2a555555 
	55888888 88555555 55a222a2 22555555 55888888 88555555 552a2a2a 2a555555 
	55888888 88555555 5522a222 a2555555 55888888 88555555 552a2a2a 2a555555 
	55888888 88555555 55a222a2 22555555 55888888 88555555 552a2a2a 2a555555 }
\shade\path(912,1212)(1512,612)(1812,912)
	(2112,612)(2712,1212)(3012,912)
	(3312,1212)(2712,1812)(2112,1212)
	(1812,1512)(1512,1212)(1212,1512)(912,1212)
\path(912,1212)(1512,612)(1812,912)
	(2112,612)(2712,1212)(3012,912)
	(3312,1212)(2712,1812)(2112,1212)
	(1812,1512)(1512,1212)(1212,1512)(912,1212)
\end{picture}}}\ar@{-}[dl]\\
&{\mbox{\begin{picture}(4224,2139)(0,-10)
\path(12,2112)(2112,12)(4212,2112)
\path(612,1512)(912,1812)(1212,1512)
	(912,1212)(612,1512)
\path(3612,1512)(3312,1812)(3012,1512)
	(3312,1212)(3612,1512)
\texture{55888888 88555555 5522a222 a2555555 55888888 88555555 552a2a2a 2a555555 
	55888888 88555555 55a222a2 22555555 55888888 88555555 552a2a2a 2a555555 
	55888888 88555555 5522a222 a2555555 55888888 88555555 552a2a2a 2a555555 
	55888888 88555555 55a222a2 22555555 55888888 88555555 552a2a2a 2a555555 }
\shade\path(912,1212)(1512,612)(1812,912)
	(2112,612)(2412,912)(2712,612)
	(3312,1212)(3012,1512)(2712,1212)
	(2412,1512)(2112,1212)(1812,1512)
	(1512,1212)(1212,1512)(912,1212)
\path(912,1212)(1512,612)(1812,912)
	(2112,612)(2412,912)(2712,612)
	(3312,1212)(3012,1512)(2712,1212)
	(2412,1512)(2112,1212)(1812,1512)
	(1512,1212)(1212,1512)(912,1212)
\end{picture}}}}}
\end{center}
The second visualization shows the graphs on $\{0,1,2,3,4,5\}$ that encode
the outer rims.
\setlength{\unitlength}{0.00014in}
\begin{center}
\mbox{\xymatrix@R=10mm@C=16mm{&
{\mbox{\begin{picture}(5424,6275)(0,-10)
\thicklines
\put(2412,3130){\makebox(0,0){\tiny$0$}}
\put(1512,3650){\makebox(0,0){\tiny$2$}}
\put(1512,2610){\makebox(0,0){\tiny$4$}}
\put(5412,3130){\makebox(0,0){\tiny$3$}}
\put(12,6248){\makebox(0,0){\tiny$5$}}
\put(12,12){\makebox(0,0){\tiny$1$}}
\path(3012,3130)(4812,3130)
\path(312,532)(4812,3130)
\path(312,532)(1212,2091)
\path(1212,4169)(1212,2091)
\path(1212,4169)(312,5728)
\path(4812,3130)(312,5728)
\end{picture}}}\ar@{-}[dl]\ar@{-}[d]\ar@{-}[dr]\\
{\mbox{\begin{picture}(5424,6275)(0,-10)
\thicklines
\put(2412,3130){\makebox(0,0){\tiny$0$}}
\put(1512,3650){\makebox(0,0){\tiny$2$}}
\put(1512,2610){\makebox(0,0){\tiny$4$}}
\put(5412,3130){\makebox(0,0){\tiny$3$}}
\put(12,6248){\makebox(0,0){\tiny$5$}}
\put(12,12){\makebox(0,0){\tiny$1$}}
\path(3012,3130)(4812,3130)
\path(312,532)(4812,3130)
\path(312,532)(1212,2091)
\path(1212,4169)(1212,2091)
\path(1212,4169)(312,5728)
\path(3012,3130)(1212,4169)
\end{picture}}}\ar@{-}[d]\ar@{-}[dr]&
{\mbox{\begin{picture}(5424,6275)(0,-10)
\thicklines
\put(2412,3130){\makebox(0,0){\tiny$0$}}
\put(1512,3650){\makebox(0,0){\tiny$2$}}
\put(1512,2610){\makebox(0,0){\tiny$4$}}
\put(5412,3130){\makebox(0,0){\tiny$3$}}
\put(12,6248){\makebox(0,0){\tiny$5$}}
\put(12,12){\makebox(0,0){\tiny$1$}}
\path(3012,3130)(4812,3130)
\path(312,532)(4812,3130)
\path(312,532)(1212,2091)
\path(312,532)(312,5728)
\path(1212,4169)(312,5728)
\path(4812,3130)(312,5728)
\end{picture}}}\ar@{-}[dl]\ar@{-}[dr]&
{\mbox{\begin{picture}(5424,6275)(0,-10)
\thicklines
\put(2412,3130){\makebox(0,0){\tiny$0$}}
\put(1512,3650){\makebox(0,0){\tiny$2$}}
\put(1512,2610){\makebox(0,0){\tiny$4$}}
\put(5412,3130){\makebox(0,0){\tiny$3$}}
\put(12,6248){\makebox(0,0){\tiny$5$}}
\put(12,12){\makebox(0,0){\tiny$1$}}
\path(3012,3130)(4812,3130)
\path(3012,3130)(1212,2091)
\path(312,532)(1212,2091)
\path(1212,4169)(1212,2091)
\path(1212,4169)(312,5728)
\path(4812,3130)(312,5728)
\end{picture}}}\ar@{-}[dl]\ar@{-}[d]\\
{\mbox{\begin{picture}(5424,6275)(0,-10)
\thicklines
\put(2412,3130){\makebox(0,0){\tiny$0$}}
\put(1512,3650){\makebox(0,0){\tiny$2$}}
\put(1512,2610){\makebox(0,0){\tiny$4$}}
\put(5412,3130){\makebox(0,0){\tiny$3$}}
\put(12,6248){\makebox(0,0){\tiny$5$}}
\put(12,12){\makebox(0,0){\tiny$1$}}
\path(3012,3130)(4812,3130)
\path(312,532)(4812,3130)
\path(312,532)(1212,2091)
\path(312,532)(312,5728)
\path(1212,4169)(312,5728)
\path(3012,3130)(1212,4169)
\end{picture}}}\ar@{-}[dr]&
{\mbox{\begin{picture}(5424,6275)(0,-10)
\thicklines
\put(2412,3130){\makebox(0,0){\tiny$0$}}
\put(1512,3650){\makebox(0,0){\tiny$2$}}
\put(1512,2610){\makebox(0,0){\tiny$4$}}
\put(5412,3130){\makebox(0,0){\tiny$3$}}
\put(12,6248){\makebox(0,0){\tiny$5$}}
\put(12,12){\makebox(0,0){\tiny$1$}}
\path(3012,3130)(4812,3130)
\path(3012,3130)(1212,2091)
\path(312,532)(1212,2091)
\path(1212,4169)(1212,2091)
\path(1212,4169)(312,5728)
\path(3012,3130)(1212,4169)
\end{picture}}}\ar@{-}[d]&
{\mbox{\begin{picture}(5424,6275)(0,-10)
\thicklines
\put(2412,3130){\makebox(0,0){\tiny$0$}}
\put(1512,3650){\makebox(0,0){\tiny$2$}}
\put(1512,2610){\makebox(0,0){\tiny$4$}}
\put(5412,3130){\makebox(0,0){\tiny$3$}}
\put(12,6248){\makebox(0,0){\tiny$5$}}
\put(12,12){\makebox(0,0){\tiny$1$}}
\path(3012,3130)(4812,3130)
\path(3012,3130)(1212,2091)
\path(312,532)(1212,2091)
\path(312,532)(312,5728)
\path(1212,4169)(312,5728)
\path(4812,3130)(312,5728)
\end{picture}}}\ar@{-}[dl]\\
&{\mbox{\begin{picture}(5424,6275)(0,-10)
\thicklines
\put(2412,3130){\makebox(0,0){\tiny$0$}}
\put(1512,3650){\makebox(0,0){\tiny$2$}}
\put(1512,2610){\makebox(0,0){\tiny$4$}}
\put(5412,3130){\makebox(0,0){\tiny$3$}}
\put(12,6248){\makebox(0,0){\tiny$5$}}
\put(12,12){\makebox(0,0){\tiny$1$}}
\path(3012,3130)(4812,3130)
\path(3012,3130)(1212,2091)
\path(312,532)(1212,2091)
\path(312,532)(312,5728)
\path(1212,4169)(312,5728)
\path(3012,3130)(1212,4169)
\end{picture}}}}}
\end{center}
\end{example}

\subsubsection*{Projection of the $1$-skeleton to the plane}
The points $j=0,\dots,N-1\in X$ correspond to the rectangular
partitions $R_j=(j^{N-j})$, where $(0^N)$ is synonymous with the empty
partition. Recall that
\begin{equation}\label{XembedEX}
f_{R_j}(k)=|k-j|(N-|k-j|).
\end{equation}
The $N$\/$\times$\/$N$ symmetric circulant matrix
$C_N:=\bigl(f_{R_j}(k)\bigr)_{j,k\in X}$
is invertible (for $N\geqslant2$); in fact, its determinant is
\begin{align*}
\det C_N
&=\prod_{j=0}^{N-1}
\Biggl(\sum_{k=0}^{N-1}k(N-k)\exp\Bigl(\frac{2\pi ijk}{N}\Bigr)\Biggr)\\
&=\frac{N(N^2-1)}{6}\prod_{j=1}^{N-1}
\frac{2N\exp\Bigl(\dfrac{2\pi ij}{N}\Bigr)}{\left(1-\exp\Bigl(
\dfrac{2\pi ij}{N}\Bigr)\right)^{\!2\mathstrut}}
=(-1)^{N+1}\dfrac{(2N)^{N-2} (N^2-1)}{3}\,,
\end{align*}
where the first equality follows by writing $\det C_N$ as the product of the
eigenvalues of $C_N$ (by using the obvious eigenvectors) and the rest is
computation (given without details).

We use the matrix
$$P:=\begin{pmatrix}
\phantom{+}0&1\\-1&0\end{pmatrix}
\begin{pmatrix}
1 & \cos\frac{2\pi}{N} & \ldots & \cos\frac{2\pi j}{N} & \ldots &
\cos\frac{2\pi (N-1)}{N} \\
0 & \sin\frac{2\pi}{N} & \ldots & \sin\frac{2\pi j}{N} & \ldots &
\sin\frac{2\pi (N-1)}{N}
\end{pmatrix}C_N^{-1}$$
to project $\mathbb R^X$ to the plane, so that $Pf_{R_0},\dots,Pf_{R_{N-1}}$
are the vertices of a regular $N$-gon. The following picture shows the
image of the $1$-skeleton of $\InjHull(X)$ for $N=9$.

\begin{center}
\includegraphics[width=10cm]{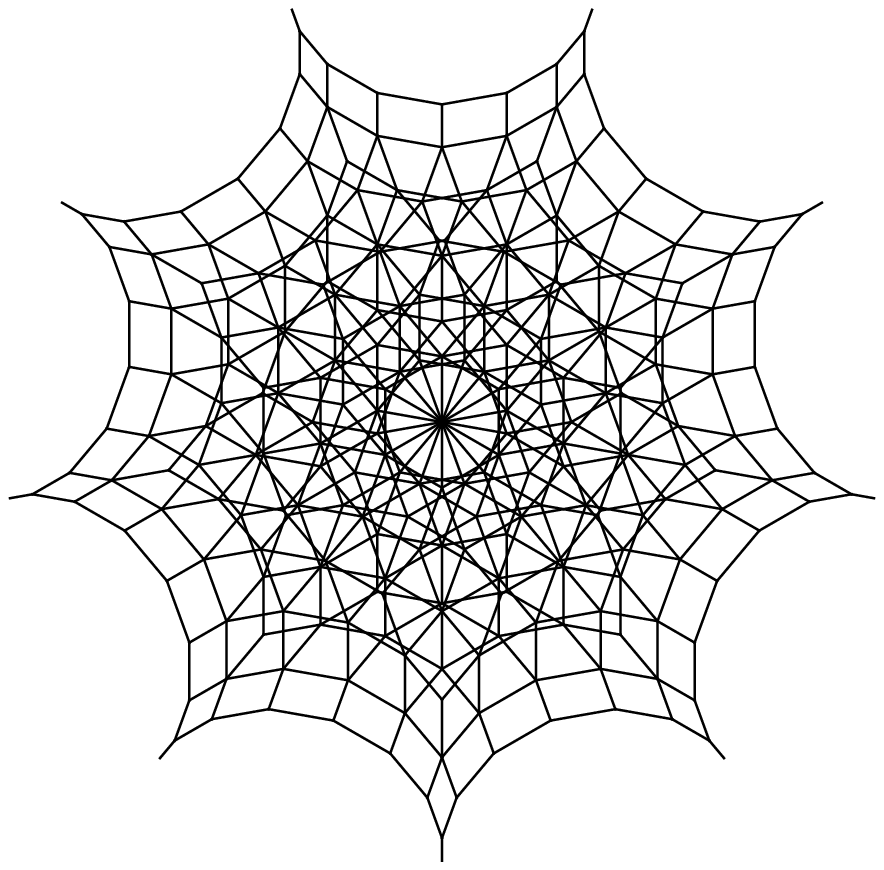}
\end{center}
We have $Pf_\lambda=(0,0)$ for the fixed point $\lambda=(4,3,2,1)$ under
the action of the cyclic group generated by
$\tau$ (see Theorem~\ref{zerocells}), but also for
$\lambda\in\bigl\{(3,3,3),(5,2,2,2),(4,4,1,1,1)\bigr\}$, which is an
orbit of size three, we have $Pf_\lambda=(0,0)$. All the other orbits
have size nine.

\section{A continuous version}\label{continuous}
Let us consider as a continuous version of integer partitions
certain functions
\begin{align}\label{continuouspartition}
\Lambda:[u,v]&\longrightarrow\mathbb R_{\geqslant0}\\\notag
u&\longmapsto -u\\\notag
v&\longmapsto v
\end{align}
with $-1\leqslant u\leqslant0\leqslant v\leqslant1$ such that the set
$\bigl\{(t,s)\bigm|t\in[u,v],\ |t|\leqslant s\leqslant\Lambda(t)\bigr\}
\subseteq\mathbb R^2$
\setlength{\unitlength}{0.00045in}
\begin{center}
\begin{picture}(6168,3480)(0,-10)
\path(15,3315)(3015,315)(6015,3315)
\path(15,315)(6015,315)
\path(1815,315)(1815,1515)
\path(4440,315)(4440,1740)
\thicklines
\path(1815,1515)(1816,1514)(1819,1512)
	(1823,1508)(1830,1502)(1838,1496)
	(1849,1489)(1863,1482)(1880,1473)
	(1902,1463)(1930,1452)(1965,1440)
	(1993,1431)(2020,1422)(2043,1415)
	(2060,1408)(2072,1403)(2080,1398)
	(2086,1394)(2090,1390)(2095,1386)
	(2101,1383)(2112,1379)(2128,1376)
	(2150,1372)(2181,1369)(2220,1366)
	(2265,1365)(2303,1366)(2339,1367)
	(2372,1370)(2400,1373)(2423,1376)
	(2443,1378)(2458,1381)(2470,1384)
	(2481,1387)(2490,1390)(2499,1393)
	(2510,1396)(2522,1400)(2537,1404)
	(2557,1408)(2580,1413)(2608,1419)
	(2641,1425)(2677,1433)(2715,1440)
	(2760,1449)(2798,1458)(2829,1465)
	(2851,1470)(2866,1475)(2875,1478)
	(2880,1481)(2884,1484)(2887,1486)
	(2891,1489)(2899,1492)(2911,1495)
	(2929,1500)(2953,1505)(2982,1510)
	(3015,1515)(3051,1520)(3080,1523)
	(3102,1525)(3116,1526)(3123,1527)
	(3127,1527)(3128,1527)(3132,1527)
	(3139,1526)(3153,1525)(3175,1523)
	(3204,1520)(3240,1515)(3273,1510)
	(3303,1504)(3327,1498)(3347,1492)
	(3361,1487)(3371,1481)(3378,1476)
	(3384,1471)(3390,1466)(3397,1461)
	(3407,1457)(3423,1452)(3443,1448)
	(3470,1444)(3503,1441)(3540,1440)
	(3577,1441)(3609,1443)(3635,1445)
	(3653,1448)(3665,1449)(3672,1451)
	(3676,1452)(3679,1454)
	(3684,1456)(3692,1459)(3705,1465)
	(3727,1472)(3757,1483)(3795,1498)
	(3840,1515)(3874,1529)(3908,1543)
	(3939,1557)(3966,1570)(3991,1582)
	(4011,1593)(4029,1603)(4043,1611)
	(4056,1619)(4067,1627)(4078,1634)
	(4088,1641)(4099,1648)(4111,1656)
	(4125,1664)(4141,1673)(4160,1683)
	(4182,1694)(4206,1705)(4233,1717)
	(4262,1729)(4290,1740)(4312,1748)
	(4332,1755)(4350,1760)(4365,1764)
	(4378,1766)(4389,1768)(4399,1769)
	(4407,1769)(4440,1740)
\thinlines
\texture{55888888 88555555 5522a222 a2555555 55888888 88555555 552a2a2a 2a555555 
	55888888 88555555 55a222a2 22555555 55888888 88555555 552a2a2a 2a555555 
	55888888 88555555 5522a222 a2555555 55888888 88555555 552a2a2a 2a555555 
	55888888 88555555 55a222a2 22555555 55888888 88555555 552a2a2a 2a555555 }
\shade\path(3015,315)(1815,1515)(1816,1514)(1819,1512)
	(1823,1508)(1830,1502)(1838,1496)
	(1849,1489)(1863,1482)(1880,1473)
	(1902,1463)(1930,1452)(1965,1440)
	(1993,1431)(2020,1422)(2043,1415)
	(2060,1408)(2072,1403)(2080,1398)
	(2086,1394)(2090,1390)(2095,1386)
	(2101,1383)(2112,1379)(2128,1376)
	(2150,1372)(2181,1369)(2220,1366)
	(2265,1365)(2303,1366)(2339,1367)
	(2372,1370)(2400,1373)(2423,1376)
	(2443,1378)(2458,1381)(2470,1384)
	(2481,1387)(2490,1390)(2499,1393)
	(2510,1396)(2522,1400)(2537,1404)
	(2557,1408)(2580,1413)(2608,1419)
	(2641,1425)(2677,1433)(2715,1440)
	(2760,1449)(2798,1458)(2829,1465)
	(2851,1470)(2866,1475)(2875,1478)
	(2880,1481)(2884,1484)(2887,1486)
	(2891,1489)(2899,1492)(2911,1495)
	(2929,1500)(2953,1505)(2982,1510)
	(3015,1515)(3051,1520)(3080,1523)
	(3102,1525)(3116,1526)(3123,1527)
	(3127,1527)(3128,1527)(3132,1527)
	(3139,1526)(3153,1525)(3175,1523)
	(3204,1520)(3240,1515)(3273,1510)
	(3303,1504)(3327,1498)(3347,1492)
	(3361,1487)(3371,1481)(3378,1476)
	(3384,1471)(3390,1466)(3397,1461)
	(3407,1457)(3423,1452)(3443,1448)
	(3470,1444)(3503,1441)(3540,1440)
	(3577,1441)(3609,1443)(3635,1445)
	(3653,1448)(3665,1449)(3672,1451)
	(3676,1452)(3679,1454)
	(3684,1456)(3692,1459)(3705,1465)
	(3727,1472)(3757,1483)(3795,1498)
	(3840,1515)(3874,1529)(3908,1543)
	(3939,1557)(3966,1570)(3991,1582)
	(4011,1593)(4029,1603)(4043,1611)
	(4056,1619)(4067,1627)(4078,1634)
	(4088,1641)(4099,1648)(4111,1656)
	(4125,1664)(4141,1673)(4160,1683)
	(4182,1694)(4206,1705)(4233,1717)
	(4262,1729)(4290,1740)(4312,1748)
	(4332,1755)(4350,1760)(4365,1764)
	(4378,1766)(4389,1768)(4399,1769)
	(4407,1769)(4440,1740)(3015,315)
\path(3015,315)(3015,3315)
\multiput(15,270)(3000,0){3}{\line(0,1){90}}
\put(2885,3315){\makebox(0,0)[r]{\footnotesize$1$}}
\put(2970,3315){\line(1,0){90}}
\put(15,-15){\makebox(0,0)[b]{\footnotesize$-1$}}
\put(6015,-15){\makebox(0,0)[b]{\footnotesize$1$}}
\put(3015,-15){\makebox(0,0)[b]{\footnotesize$0$}}
\put(1815,-15){\makebox(0,0)[b]{\footnotesize$u$}}
\put(4440,-15){\makebox(0,0)[b]{\footnotesize$v$}}
\end{picture}
\end{center}
has an evident meaning of ``continuous Young diagram'' of the
``continuous partition''~$\Lambda$.
That the parts of an integer partition are ordered
decreasingly translates into the requirement that $\Lambda$
must be $1$-Lipschitz. We further require that $\Lambda$ be
defined on an interval $[u,v]$ of length at most $1$ (which
replaces the bound on the hook lengths for integer partitions
in $\mathbb Y_N$). There is then a unique extension of $\Lambda$
to a $1$-Lipschitz function $\Lambda:[u,u+1]\to\mathbb R_{\geqslant0}$
with $\Lambda(u+1)=u+1$, namely, $\Lambda(t)=t$ for $t\in[v,u+1]$.
Note that $\Lambda(t)\leqslant1$ for all $t\in[u,u+1]$ because
for $u<t<u+1$ the $1$-Lipschitz condition gives
$$\frac{\Lambda(t)+u}{t-u}=\frac{\Lambda(t)-\Lambda(u)}{t-u}\leqslant1
\mbox{\quad and\quad}
\frac{\Lambda(t)-u-1}{u+1-t}=\frac{\Lambda(t)-\Lambda(u+1)}{u+1-t}\leqslant1$$
and hence
$$2\Lambda(t)-1=\bigl(\Lambda(t)+u\bigr)+\bigl(\Lambda(t)-u-1\bigr)\leqslant
(t-u)+(u+1-t)=1.$$
The shifted antiperiodic extension defined by $\Lambda(t+1)=1-\Lambda(t)$
finally extends $\Lambda$ to the real line $\Lambda:\mathbb R\to[0,1]$,
so that the graph of $\Lambda$ becomes the preimage of a loop in the
M\"obius strip under the universal covering projection.
We define
$$\mathbb Y_\infty:=\bigl\{\Lambda:\mathbb R\to[0,1]\bigm|
\mbox{$\Lambda$ is $1$-Lipschitz and $\Lambda(t+1)=1-\Lambda(t)$ for all
$t$}\bigr\}.$$
For each $\Lambda\in\mathbb Y_\infty$ its graph intersects the line segment
with endpoints $(-1,1)$ and $(0,0)$ in at least one (and generically 
exactly one) point.
If $(u,-u)$ is a point in the intersection, then the
restriction $\Lambda|_{[u,u+1]}$ defines a ``continuous partition''
as in (\ref{continuouspartition}) with $v=u+1$.

\setlength{\unitlength}{0.00045in}
\begin{center}
\begin{picture}(13824,3680)(0,-10)
\path(1812,2115)(3012,3315)(4318,2009)
	(6012,315)(7318,1621)
\path(7318,1621)(9012,3315)
\path(9012,3315)(10812,1515)
\thicklines
\path(9012,1515)(9014,1515)(9019,1515)
	(9028,1516)(9042,1517)(9060,1518)
	(9084,1520)(9112,1522)(9144,1525)
	(9179,1528)(9215,1532)(9253,1536)
	(9291,1541)(9328,1546)(9364,1552)
	(9400,1558)(9434,1565)(9467,1573)
	(9498,1581)(9529,1590)(9560,1600)
	(9590,1611)(9620,1623)(9650,1636)
	(9681,1650)(9712,1665)(9737,1677)
	(9762,1691)(9787,1705)(9814,1719)
	(9841,1734)(9869,1750)(9898,1767)
	(9927,1784)(9958,1802)(9989,1821)
	(10021,1839)(10053,1859)(10087,1878)
	(10121,1898)(10155,1918)(10190,1938)
	(10226,1958)(10262,1978)(10297,1998)
	(10333,2017)(10369,2036)(10405,2054)
	(10441,2072)(10476,2089)(10511,2105)
	(10546,2121)(10580,2136)(10615,2150)
	(10648,2163)(10681,2175)(10714,2187)
	(10747,2197)(10780,2206)(10812,2215)
	(10844,2223)(10877,2229)(10910,2235)
	(10943,2240)(10976,2243)(11011,2246)
	(11047,2247)(11083,2248)(11122,2248)
	(11161,2246)(11203,2244)(11246,2241)
	(11291,2237)(11338,2232)(11387,2226)
	(11437,2219)(11489,2211)(11506,2209)
\path(2506,1421)(2542,1427)
	(2595,1436)(2648,1445)(2699,1454)
	(2750,1463)(2797,1472)(2841,1480)
	(2880,1488)(2915,1495)(2944,1501)
	(2967,1505)(2985,1509)(2998,1512)
	(3006,1514)(3010,1515)(3012,1515)
\path(3012,1515)(3014,1515)(3019,1515)
	(3028,1516)(3042,1517)(3060,1518)
	(3084,1520)(3112,1522)(3144,1525)
	(3179,1528)(3215,1532)(3253,1536)
	(3291,1541)(3328,1546)(3364,1552)
	(3400,1558)(3434,1565)(3467,1573)
	(3498,1581)(3529,1590)(3560,1600)
	(3590,1611)(3620,1623)(3650,1636)
	(3681,1650)(3712,1665)(3737,1677)
	(3762,1691)(3787,1705)(3814,1719)
	(3841,1734)(3869,1750)(3898,1767)
	(3927,1784)(3958,1802)(3989,1821)
	(4021,1839)(4053,1859)(4087,1878)
	(4121,1898)(4155,1918)(4190,1938)
	(4226,1958)(4262,1978)(4297,1998)(4318,2009)
	(4333,2017)(4369,2036)(4405,2054)
	(4441,2072)(4476,2089)(4511,2105)
	(4546,2121)(4580,2136)(4615,2150)
	(4648,2163)(4681,2175)(4714,2187)
	(4747,2197)(4780,2206)(4812,2215)
	(4844,2223)(4877,2229)(4910,2235)
	(4943,2240)(4976,2243)(5011,2246)
	(5047,2247)(5083,2248)(5122,2248)
	(5161,2246)(5203,2244)(5246,2241)
	(5291,2237)(5338,2232)(5387,2226)
	(5437,2219)(5489,2211)(5542,2203)
	(5595,2194)(5648,2185)(5699,2176)
	(5750,2167)(5797,2158)(5841,2150)
	(5880,2142)(5915,2135)(5944,2129)
	(5967,2125)(5985,2121)(5998,2118)
	(6006,2116)(6010,2115)(6012,2115)
\thinlines
\texture{55888888 88555555 5522a222 a2555555 55888888 88555555 552a2a2a 2a555555 
	55888888 88555555 55a222a2 22555555 55888888 88555555 552a2a2a 2a555555 
	55888888 88555555 5522a222 a2555555 55888888 88555555 552a2a2a 2a555555 
	55888888 88555555 55a222a2 22555555 55888888 88555555 552a2a2a 2a555555 }
\shade\path(6012,315)(4318,2009)
	(4333,2017)(4369,2036)(4405,2054)
	(4441,2072)(4476,2089)(4511,2105)
	(4546,2121)(4580,2136)(4615,2150)
	(4648,2163)(4681,2175)(4714,2187)
	(4747,2197)(4780,2206)(4812,2215)
	(4844,2223)(4877,2229)(4910,2235)
	(4943,2240)(4976,2243)(5011,2246)
	(5047,2247)(5083,2248)(5122,2248)
	(5161,2246)(5203,2244)(5246,2241)
	(5291,2237)(5338,2232)(5387,2226)
	(5437,2219)(5489,2211)(5542,2203)
	(5595,2194)(5648,2185)(5699,2176)
	(5750,2167)(5797,2158)(5841,2150)
	(5880,2142)(5915,2135)(5944,2129)
	(5967,2125)(5985,2121)(5998,2118)
	(6006,2116)(6010,2115)(6012,2115)
        (6014,2115)(6019,2115)
	(6028,2114)(6042,2113)(6060,2112)
	(6084,2110)(6112,2108)(6144,2105)
	(6179,2102)(6215,2098)(6253,2094)
	(6291,2089)(6328,2084)(6364,2078)
	(6400,2072)(6434,2065)(6467,2057)
	(6498,2049)(6529,2040)(6560,2030)
	(6590,2019)(6620,2007)(6650,1994)
	(6681,1980)(6712,1965)(6737,1953)
	(6762,1939)(6787,1925)(6814,1911)
	(6841,1896)(6869,1880)(6898,1863)
	(6927,1846)(6958,1828)(6989,1809)
	(7021,1791)(7053,1771)(7087,1752)
	(7121,1732)(7155,1712)(7190,1692)
	(7226,1672)(7262,1652)(7297,1632)(7318,1621)
        (6012,315)
\path(6012,315)(4318,2009)
	(4333,2017)(4369,2036)(4405,2054)
	(4441,2072)(4476,2089)(4511,2105)
	(4546,2121)(4580,2136)(4615,2150)
	(4648,2163)(4681,2175)(4714,2187)
	(4747,2197)(4780,2206)(4812,2215)
	(4844,2223)(4877,2229)(4910,2235)
	(4943,2240)(4976,2243)(5011,2246)
	(5047,2247)(5083,2248)(5122,2248)
	(5161,2246)(5203,2244)(5246,2241)
	(5291,2237)(5338,2232)(5387,2226)
	(5437,2219)(5489,2211)(5542,2203)
	(5595,2194)(5648,2185)(5699,2176)
	(5750,2167)(5797,2158)(5841,2150)
	(5880,2142)(5915,2135)(5944,2129)
	(5967,2125)(5985,2121)(5998,2118)
	(6006,2116)(6010,2115)(6012,2115)
        (6014,2115)(6019,2115)
	(6028,2114)(6042,2113)(6060,2112)
	(6084,2110)(6112,2108)(6144,2105)
	(6179,2102)(6215,2098)(6253,2094)
	(6291,2089)(6328,2084)(6364,2078)
	(6400,2072)(6434,2065)(6467,2057)
	(6498,2049)(6529,2040)(6560,2030)
	(6590,2019)(6620,2007)(6650,1994)
	(6681,1980)(6712,1965)(6737,1953)
	(6762,1939)(6787,1925)(6814,1911)
	(6841,1896)(6869,1880)(6898,1863)
	(6927,1846)(6958,1828)(6989,1809)
	(7021,1791)(7053,1771)(7087,1752)
	(7121,1732)(7155,1712)(7190,1692)
	(7226,1672)(7262,1652)(7297,1632)(7318,1621)
        (6012,315)
\thicklines
\path(6012,2115)(6014,2115)(6019,2115)
	(6028,2114)(6042,2113)(6060,2112)
	(6084,2110)(6112,2108)(6144,2105)
	(6179,2102)(6215,2098)(6253,2094)
	(6291,2089)(6328,2084)(6364,2078)
	(6400,2072)(6434,2065)(6467,2057)
	(6498,2049)(6529,2040)(6560,2030)
	(6590,2019)(6620,2007)(6650,1994)
	(6681,1980)(6712,1965)(6737,1953)
	(6762,1939)(6787,1925)(6814,1911)
	(6841,1896)(6869,1880)(6898,1863)
	(6927,1846)(6958,1828)(6989,1809)
	(7021,1791)(7053,1771)(7087,1752)
	(7121,1732)(7155,1712)(7190,1692)
	(7226,1672)(7262,1652)(7297,1632)(7318,1621)
	(7333,1613)(7369,1594)(7405,1576)
	(7441,1558)(7476,1541)(7511,1525)
	(7546,1509)(7580,1494)(7615,1480)
	(7648,1467)(7681,1455)(7714,1443)
	(7747,1433)(7780,1424)(7812,1415)
	(7844,1407)(7877,1401)(7910,1395)
	(7943,1390)(7976,1387)(8011,1384)
	(8047,1383)(8083,1382)(8122,1382)
	(8161,1384)(8203,1386)(8246,1389)
	(8291,1393)(8338,1398)(8387,1404)
	(8437,1411)(8489,1419)(8542,1427)
	(8595,1436)(8648,1445)(8699,1454)
	(8750,1463)(8797,1472)(8841,1480)
	(8880,1488)(8915,1495)(8944,1501)
	(8967,1505)(8985,1509)(8998,1512)
	(9006,1514)(9010,1515)(9012,1515)
\thinlines
\path(6012,2115)(6012,315)%  t
\path(4318,2009)(4318,315)%  a
\path(7318,1621)(7318,315)%  b
\path(612,3315)(12612,3315)
\path(2506,315)
     (9612,315)
\path(2506,315)(2506,1421)% u
\path(5506,315)(5506,2209)% u+1
\path(612,3315)(3612,315)(6612,3315)
	(9612,315)(12612,3315)
\put(6612,3385){\makebox(0,0)[b]{\footnotesize$(1,1)$}}
\put(12612,3385){\makebox(0,0)[b]{\footnotesize$(3,1)$}}
\put(612,3385){\makebox(0,0)[b]{\footnotesize$(-1,1)\phantom{\ }$}}
\put(3612,215){\makebox(0,0)[t]{\footnotesize$(0,0)$}}
\put(9612,215){\makebox(0,0)[t]{\footnotesize$(2,0)$}}
\put(6012,15){\makebox(0,0)[b]{\footnotesize$t$}}
\put(5865,2215){\makebox(0,0)[lb]{\footnotesize$\Lambda(t)$}}
\put(4318,15){\makebox(0,0)[b]{\footnotesize$a$}}
\put(7318,15){\makebox(0,0)[b]{\footnotesize$b=a\raisebox{0pt}[0pt][0pt]{$+$}1$}}
\put(5506,15){\makebox(0,0)[b]{\footnotesize$u\raisebox{0pt}[0pt][0pt]{$+$}1$}}
\put(2506,15){\makebox(0,0)[b]{\footnotesize$u$}}
\put(4378,2309){\makebox(0,0)[b]{\footnotesize$\Lambda(a)$}}
\put(7300,1921){\makebox(0,0)[b]{\footnotesize$\Lambda(b)$}}
\put(2400,1320){\makebox(0,0)[r]{\footnotesize$\Lambda(u)=-u$}}
\dashline[24]{80}(4318,2009)(5624,3315)(7318,1621)(8624,315)(10318,2009)
\put(8624,15){\makebox(0,0)[b]{\footnotesize$s$}}
\put(5623,3385){\makebox(0,0)[b]{\footnotesize$(s\!-\!1,1)\phantom{+0}$}}
\put(11606,2209){\makebox(0,0)[l]{\footnotesize{graph of $\Lambda$}}}
\end{picture}
\end{center}

In analogy with the assignment $\mathbb Y_N\ni\lambda\mapsto f_\lambda$ we
define $\mathbb Y_\infty\ni\Lambda\mapsto F_\Lambda$ by the ($2$-periodic)
area function
$$F_\Lambda(t)=\int_a^b \Lambda(s)\,ds
-\frac12\bigl(\Lambda(a)^2+\Lambda(b)^2\bigr)$$
where $b=a+1$ and $t=a+\Lambda(a)$ and hence $t=b-\Lambda(b)$;
moreover, $s=b+\Lambda(b)$.

\mbox{}\\\indent
For $u\leqslant a\leqslant u+1$ and hence
$0\leqslant a+\Lambda(a)\leqslant2+2u\leqslant a+1+\Lambda(a+1)\leqslant2$
we have
\begin{align}\notag
F_\Lambda(t)&=F_\Lambda\bigl(a+\Lambda(a)\bigr)=\int_a^{u+1}\Lambda(s)\,ds+
\int_{u+1}^{a+1}\Lambda(s)\,ds
-\frac12\bigl(\Lambda(a)^2+\Lambda(a+1)^2\bigr)\\\notag
&=\int_a^{u+1}\Lambda(s)\,ds+\int_u^a\bigl(1-\Lambda(s)\bigr)\,ds
-\frac12\bigl(\Lambda(a)^2+(1-\Lambda(a))^2\bigr)\\\label{Fa}
&=a-u-\int_u^a\Lambda(s)\,ds+\int_a^{u+1}\Lambda(s)\,ds
+\Lambda(a)\bigl(1-\Lambda(a)\bigr)-\frac12\,,\\[1mm]\notag
F_\Lambda(s)&=F_\Lambda\bigl(b+\Lambda(b)\bigr)
=F_\Lambda\bigl(a+1+\Lambda(a+1)\bigr)=
F_\Lambda\bigl(2+a-\Lambda(a)\bigr)\phantom{\int_0^0}\\\notag
&=\int_{a+1}^{u+2}\Lambda(s)\,ds+\int_{u+2}^{a+2}\Lambda(s)\,ds
-\frac12\bigl(\Lambda(a+1)^2+\Lambda(a+2)^2\bigr)\\\notag
&=\int_a^{u+1}\bigl(1-\Lambda(s)\bigr)\,ds+\int_u^a\Lambda(s)\,ds
-\frac12\bigl((1-\Lambda(a))^2+\Lambda(a)^2\bigr)\\\label{Fb}
&=u+1-a+\int_u^a\Lambda(s)\,ds-\int_a^{u+1}\Lambda(s)\,ds
+\Lambda(a)\bigl(1-\Lambda(a)\bigr)-\frac12\,.
\end{align}
The equality
\begin{equation}\label{complement}
F_\Lambda\bigl(a+\Lambda(a)\bigr)+F_\Lambda\bigl(b+\Lambda(b)\bigr)
=2\Lambda(a)\Lambda(b)\qquad\mbox{(for $b=a+1$)}
\end{equation}
follows directly from the previous picture by considering the rectangle
with corners at $(a,\Lambda(a))$, $(t,0)$, $(b,\Lambda(b))$, $(s-1,1)$;
alternatively, (\ref{complement}) follows by adding (\ref{Fa}) and (\ref{Fb}).

For $r\in[0,1]$ let $R_r\in\mathbb Y_\infty$ be the function that satisfies
$$R_r(t)=\begin{cases}
t+2-2r&\mbox{if $r-1\leqslant t\leqslant2r-1$},\\
-t+2r&\mbox{if $2r-1\leqslant t\leqslant r$}.
\end{cases}$$
\begin{center}
\begin{picture}(9134,3352)(0,-10)
\path(912,315)(9012,315)
\path(912,315)(912,2415)
\path(3912,315)(3912,1215)
\path(1812,315)(1815,3315)
\thicklines
\path(912,2415)(1812,3315)(4812,315)(6912,2415)
\thinlines
\path(6912,315)(6912,2415)
\path(12,3315)(3012,315)(6012,3315)(9012,315)
\path(12,3315)(6012,3315)
\put(912,15){\makebox(0,0)[t]{\footnotesize$r\!-\!1$}}
\put(3012,15){\makebox(0,0)[t]{\footnotesize$0$}}
\put(3912,15){\makebox(0,0)[t]{\footnotesize$\phantom{1}r\phantom{1}$}}
\put(4812,15){\makebox(0,0)[t]{\footnotesize$2r$}}
\put(1812,15){\makebox(0,0)[t]{\footnotesize$2r\!-\!1$}}
\put(6912,15){\makebox(0,0)[t]{\footnotesize$r+1$}}
\put(9012,15){\makebox(0,0)[t]{\footnotesize$2$}}
\put(7212,2415){\makebox(0,0)[l]{\footnotesize graph of $ R_r$}}
\end{picture}
\end{center}
Note that $R_0=R_1$. Note also that $R_r$ is the unique function in
$\mathbb Y_\infty$ that vanishes at $2r$. 
The $2$-periodic function $F_{R_r}$ is then
$$F_{R_r}(t)=\frac12(t-2r)(2r+2-t)\qquad\mbox{for $2r\leqslant t\leqslant
2r+2$}.$$
We define the metric space $X_\infty$ as the subspace
$$X_\infty:=\bigl\{R_r\bigm|r\in[0,1]\bigr\}\subseteq\mathbb Y_\infty.$$
Two alternative characterizations are
$$X_\infty=\bigl\{\Lambda\in\mathbb Y_\infty\bigm|
\exists\,t\in\mathbb R: \Lambda(t)=0\bigr\}
=\bigl\{\Lambda\in\mathbb Y_\infty\bigm|
\exists\,t\in\mathbb R: F_\Lambda(t)=0\bigr\}.$$
The distance function is
$$D(r,s):=d(R_r,R_s)=\Vert F_{R_r}-F_{R_s}\Vert_\infty=2|s-r|\bigl(1-|s-r|\bigr).$$
For $r,s\in[0,1]$ the equality
$$F_{R_r}(2s)=2|s-r|\bigl(1-|s-r|\bigr)$$
is an analogue of (\ref{XembedEX}).

\subsubsection*{A short dictionary}
\renewcommand{\arraystretch}{1.4}
$$\begin{array}{|l@{\quad}|@{\quad}l|}\hline
\mathbb Y_N\qquad\mbox{($N\in\mathbb Z_{\geqslant 2}$)}&\mathbb Y_\infty\\\hline
\lambda\in\mathbb Y_N&\Lambda:[u,u+1]\to[0,1]\mbox{ $1$-Lipschitz}\\[-1mm]
&\mbox{with $\Lambda(u)=-u$ and $\Lambda(u+1)=u+1$}\\[-1mm]
&\mbox{extended by $\Lambda(t+1)=1-\Lambda(t)$ to all $t\in\mathbb R$}\\\hline
(0,0),\,(0,N),\,(N,N);\,(j,k)&(-1,1),\,(0,0),\,(1,1);\,
\bigl(\frac jN+\frac kN-1,\frac jN-\frac kN+1\bigr)\\\hline
(j,k)\in\mathcal L_\lambda\mbox{ (for $0\leqslant j\leqslant k\leqslant N$)}
&\bigl(a,\Lambda(a)\bigr)\mbox{ (for $u\leqslant a
\leqslant u+1$)}\\\hline
d(j,k)=|k-j|\bigl(N-|k-j|\bigr)&D(r,s)=2|s-r|\bigl(1-|s-r|\bigr)\\\hline
\forall\,(j,k)\in\mathcal L_\lambda:&\mbox{for $b=a+1$:}\\[-1mm]
f_\lambda(j)+f_\lambda(k)=d(j,k)&
F_\Lambda\bigl(a+\Lambda(a)\bigr)+F_\Lambda\bigl(b+\Lambda(b)\bigr)
=2\Lambda(a)\Lambda(b)\\[-1mm]
{}=d(0,|k-j|)&{}=D(0,\Lambda(a))=D(0,\Lambda(b))\\\hline
X_N&X_\infty\\\hline
\end{array}$$
\mbox{}\\
Of course one could elaborate on the correspondence by first putting
the representation of integer partitions into a form analogous to
the representation of ``continuous partitions''. For instance, look at
the partition $\lambda=(5,3,3,2)\in\mathbb Y_N$ with $N=9$.

\setlength{\unitlength}{0.00057in}
\begin{center}
\begin{picture}(10824,4085)(0,-25)
\path(12,3462)(312,3762)(612,3462)
	(912,3762)(1212,3462)(1512,3762)
	(1812,3462)(2112,3762)(2412,3462)
	(2712,3762)(3012,3462)(3312,3762)
	(3612,3462)(3912,3762)(4212,3462)
	(4512,3762)(4812,3462)(5112,3762)
	(5412,3462)(5712,3762)(6012,3462)
	(6312,3762)(6612,3462)(6912,3762)
	(7212,3462)(7512,3762)(7812,3462)
	(8112,3762)(8412,3462)(8712,3762)
	(9012,3462)(9312,3762)(9612,3462)
	(9912,3762)(10212,3462)(10512,3762)
	(10812,3462)(7812,462)(7512,762)
	(7212,462)(6912,762)(6612,462)
	(6312,762)(6012,462)(5712,762)
	(5412,462)(5112,762)(4812,462)
	(4512,762)(4212,462)(3912,762)
	(3612,462)(3312,762)(3012,462)(12,3462)
\path(1512,1962)(2412,2862)(2712,2562)
	(3012,2862)(3612,2262)(4212,2862)
	(5112,1962)(5412,2262)(5712,1962)
	(6312,2562)(6912,1962)(7812,2862)
	(8112,2562)(8412,2862)(9012,2262)
	(9612,2862)(9912,2562)(9012,1662)
	(8412,2262)(8112,1962)(7812,2262)
	(6912,1362)(6312,1962)(5712,1362)
	(5412,1662)(5112,1362)(4212,2262)
	(3612,1662)(3012,2262)(2712,1962)
	(2412,2262)(1812,1662)(1512,1962)
\thicklines
\path(1812,1662)(3012,462)(4512,1962)
	(4212,2262)(3612,1662)(3012,2262)
	(2712,1962)(2412,2262)(1812,1662)
\thinlines
\path(1512,762)(10812,762)
\path(10692.000,732.000)(10812.000,762.000)(10692.000,792.000)
\path(3012,762)(3012,2412)
\path(3012,2712)(3012,3762)
\path(3042.000,3642.000)(3012.000,3762.000)(2982.000,3642.000)
\path(2112,1362)(2712,1962)
\path(2412,1062)(3312,1962)
\path(2712,762)(3612,1662)
\path(2112,1962)(3312,762)
\path(2712,1962)(3612,1062)
\path(3612,1662)(3912,1362)
\path(3912,1962)(4212,1662)
\path(1662,162)(1662,12)(4362,12)
	(4362,162)(4362,12)(7062,12)(7062,162)
\path(4212,2262)(4512,2562)
\path(6912,1962)(7212,1662)
\put(1812,162){\makebox(0,0){\tiny$-4$}}
\put(2112,162){\makebox(0,0){\tiny$-3$}}
\put(2412,162){\makebox(0,0){\tiny$-2$}}
\put(2712,162){\makebox(0,0){\tiny$-1$}}
\put(3012,162){\makebox(0,0){\tiny$0$}}
\put(3312,162){\makebox(0,0){\tiny$1$}}
\put(3612,162){\makebox(0,0){\tiny$2$}}
\put(3912,162){\makebox(0,0){\tiny$3$}}
\put(4212,162){\makebox(0,0){\tiny$4$}}
\put(4512,162){\makebox(0,0){\tiny$5$}}
\put(4812,162){\makebox(0,0){\tiny$6$}}
\put(5112,162){\makebox(0,0){\tiny$7$}}
\put(5412,162){\makebox(0,0){\tiny$8$}}
\put(5712,162){\makebox(0,0){\tiny$9$}}
\put(6012,162){\makebox(0,0){\tiny$10$}}
\put(6312,162){\makebox(0,0){\tiny$11$}}
\put(6612,162){\makebox(0,0){\tiny$12$}}
\put(6912,162){\makebox(0,0){\tiny$13$}}
\put(7212,162){\makebox(0,0){\tiny$14$}}
\put(7512,162){\makebox(0,0){\tiny$15$}}
\put(7812,162){\makebox(0,0){\tiny$16$}}
\put(1812,1962){\makebox(0,0){\footnotesize$4$}}
\put(2112,2262){\makebox(0,0){\footnotesize$5$}}
\put(2412,2562){\makebox(0,0){\footnotesize$6$}}
\put(2712,2262){\makebox(0,0){\footnotesize$5$}}
\put(3012,2562){\makebox(0,0){\footnotesize$6$}}
\put(3312,2262){\makebox(0,0){\footnotesize$5$}}
\put(3612,1962){\makebox(0,0){\footnotesize$4$}}
\put(3912,2262){\makebox(0,0){\footnotesize$5$}}
\put(4212,2562){\makebox(0,0){\footnotesize$6$}}
\put(4512,2262){\makebox(0,0){\footnotesize$5$}}
\put(4812,1962){\makebox(0,0){\footnotesize$4$}}
\put(5112,1662){\makebox(0,0){\footnotesize$3$}}
\put(5412,1962){\makebox(0,0){\footnotesize$4$}}
\put(5712,1662){\makebox(0,0){\footnotesize$3$}}
\put(6012,1962){\makebox(0,0){\footnotesize$4$}}
\put(6312,2262){\makebox(0,0){\footnotesize$5$}}
\put(6612,1962){\makebox(0,0){\footnotesize$4$}}
\put(6912,1662){\makebox(0,0){\footnotesize$3$}}
\put(7212,1962){\makebox(0,0){\footnotesize$4$}}
\put(7512,2262){\makebox(0,0){\footnotesize$5$}}
\put(7812,2562){\makebox(0,0){\footnotesize$6$}}
\put(8112,2262){\makebox(0,0){\footnotesize$5$}}
\put(8412,2562){\makebox(0,0){\footnotesize$6$}}
\put(8712,2262){\makebox(0,0){\footnotesize$5$}}
\put(9012,1962){\makebox(0,0){\footnotesize$4$}}
\put(9312,2262){\makebox(0,0){\footnotesize$5$}}
\put(9612,2562){\makebox(0,0){\footnotesize$6$}}
\put(8112,162){\makebox(0,0){\tiny$17$}}
\put(8412,162){\makebox(0,0){\tiny$18$}}
\put(8712,162){\makebox(0,0){\tiny$19$}}
\put(9012,162){\makebox(0,0){\tiny$20$}}
\put(9312,162){\makebox(0,0){\tiny$21$}}
\put(9612,162){\makebox(0,0){\tiny$22$}}
\end{picture}
\end{center}
The function $\Lambda_9:[-4,4]\cap\mathbb Z\to[0,9]\cap\mathbb Z$\/
$$\begin{array}{c|rrrrrrrrr}
t&-4&-3&-2&-1&\phantom{+}0&\phantom{+}1&\phantom{+}2&\phantom{+}3&
\phantom{+}4\\\hline
\Lambda_9(t)&4&5&6&5&6&5&4&5&6
\end{array}$$
describes the heights of the box positions in the outer rim
of $\lambda$. In general, for $\lambda\in\mathbb Y_N$ we have the
corresponding function $\Lambda_N$ defined on a set
of $N$ consecutive integers that contains $0$, and we extend it to
a function $\Lambda_N:\mathbb Z\to[0,N]\cap\mathbb Z$
by requiring the symmetry $\Lambda_N(t+N)=N-\Lambda_N(t)$.
In the example we have $\Lambda_9(-4)=4$ and $\Lambda_9(5)=5$, which is
reminiscent of the boundary condition $\Lambda(u)=-u$ and $\Lambda(u+1)=u+1$
for a ``continuous partition'' $\Lambda:[u,u+1]\to[0,1]$.
Let us also write down the formula
$$\Lambda_N(0)=2\cdot\sqrt{\#\mathstrut
(\mbox{boxes in the Durfee square of the partition $\lambda$})}.$$

\section{Some conclusions and outlook}
The realization of $\operatorname{Hasse}(\mathbb Y_N)$ as the $1$-skeleton of
the $\bigl\lfloor\frac N2\bigr\rfloor$-dimensional polyhedral complex
$\InjHull(X_N)$ embedded in $\mathbb R^N$ (in fact, in its nonnegative
orthant) that is stable under cyclic permutation of the coordinates, provides
a new geometric realization of the $N$-fold cyclic symmetry of
$\operatorname{Hasse}(\mathbb Y_N)$.

When one looks at the graphs $\operatorname{Hasse}(\mathbb Y_N)$ depicted in
\cite{Su1} for $N=5,6,7,8$, then the $1$-skeleta of five squares, one cube,
seven cubes, respectively one tesseract are eye-catching. To some extent,
the construction of $\InjHull(X_N)$ gives a meaning to those
$\bigl\lfloor\frac N2\bigr\rfloor$-cubes. In general, Theorem~\ref{allfaces}
counts all those as well as the lower-dimensional $v$-cubes.

The works of C.~Berg and M.~Zabrocki \cite{BZ} and H.~Thomas and N.~Williams
\cite{TW} generalize the cyclic symmetries and prove an instance of the cyclic
sieving phenomenon in that framework. One may ask whether the metric geometry
approach developed here for the original cyclic symmetries can be extended
in some nice way to the generalized version.

Homotopically nontrivial loops (as specified in
Definition~\ref{shortnontrivialloops}) play an important role.
An investigation of loops (self-avoiding or not, of given length and with a
prescribed winding number) in various sorts of discrete strips (also with
other than $\mathsf A_{N-1}$ type Dynkin diagram fibres) might involve
interesting combinatorics.

Let us conclude with some geometry.
The maximal coordinate sum of points in $\InjHull(X_N)$ is attained for the
$N$ points in the orbit of $f_{()}$,
\begin{align*}
\max\bigl\{\Vert f\Vert_1\bigm|f\in\InjHull(X_N)\bigr\}
&=\Vert f_{()}\Vert_1=\sum_{j=0}^{N-1}j(N-j)=\frac16(N^3-N).
\intertext{For $N$ odd there is a unique partition
$\lambda\in\mathbb Y_N$ with}
\min\bigl\{\Vert f\Vert_1\bigm|f\in\InjHull(X_N)\bigr\}
&=\Vert f_\lambda\Vert_1=\frac18(N^3-N), 
\intertext{namely the staircase partition
$\lambda=\bigl(\frac{N-1}{2},\frac{N-3}{2},\dots,1\bigr)$, such that
$f_\lambda(l)=\frac18(N^2-1)$ for all $l$.
For $N$ even there is a central $\frac N2$-dimensional cube, whose points
have minimal $1$-norm}
\min\bigl\{\Vert f\Vert_1\bigm|f\in\InjHull(X_N)\bigr\}
&=\frac18 N^3.
\end{align*}
The vertices of this $\frac N2$-cube are $f_\nu$ where $\nu$ is one of the
$2^{N/2}$ partitions that are got by removing from the staircase partition
$\bigl(\frac N2,\frac{N-2}{2},\dots,1\bigr)$ any subset of its inner corners.
The barycentre of this $\frac N2$-cube has all coordinates $\frac18 N^2$,
but no vertex of $\operatorname{Hasse}(\mathbb Y_N)$ is fixed under the
cyclic action (which is also a very particular case of the cyclic sieving
phenomenon established in \cite{TW}).

Recall that it is the supremum norm in $\mathbb R^N$ that endows
$\InjHull(X_N)$ with the metric that characterizes this space as an
injective hull. But the polyhedral complex $\InjHull(X_N)$ can also be
looked at from a more geometric point of view by using the Euclidean metric,
so that all the edges have length $\sqrt{N}$ by Theorem~\ref{zerocells}.
In the Euclidean metric the `central cube' for $N=6$ is actually a
rhombohedron with six congruent faces, namely, rhombi with diagonals of
lengths $2\sqrt{2}$ and $4$; one of the diagonals of the rhombohedron has
length $\sqrt{6}$ and the other three have length $\sqrt{22}$ $\Bigl($which is
$\sqrt{\frac{1}{12}(N^3+8N)}\Bigr|_{N=6}\Bigr)$.

\end{document}